\documentclass[11pt,a4paper]{article}
\usepackage[margin=2.6cm]{geometry}
\usepackage{amsmath,amssymb,amsthm}
\usepackage{booktabs}
\usepackage{array}
\usepackage{tikz}
\usetikzlibrary{positioning,calc}
\usepackage[colorlinks=true, linkcolor=blue!60!black, citecolor=blue!60!black, urlcolor=blue!60!black]{hyperref}
\usepackage{microtype}
\usepackage{enumitem}
\usepackage{url}
\usepackage[numbers,sort&compress]{natbib}
\usepackage{adjustbox}

\newtheorem{theorem}{Theorem}[section]
\newtheorem{lemma}[theorem]{Lemma}
\newtheorem{proposition}[theorem]{Proposition}
\newtheorem{corollary}[theorem]{Corollary}
\newtheorem{conjecture}[theorem]{Conjecture}
\newtheorem{introtheorem}{Theorem}

\theoremstyle{definition}
\newtheorem{definition}[theorem]{Definition}
\newtheorem{example}[theorem]{Example}

\theoremstyle{remark}
\newtheorem{remark}[theorem]{Remark}

\newcommand{\R}{\mathbb{R}}
\newcommand{\Z}{\mathbb{Z}}
\newcommand{\Lp}{{L^{\dagger}}}
\newcommand{\tr}{\operatorname{tr}}
\newcommand{\Kf}{\operatorname{Kf}}
\newcommand{\one}{\mathbf{1}}
\newcommand{\rres}[2]{r_{#1 #2}}
\newcommand{\Fp}{\mathcal{F}}

\title{Resistance-regular graphs need not be regular}

\author{%
	Murad ul Islam Khan\thanks{Email: \texttt{muradulislam@uoh.edu.pk}}
	\quad and \quad
	Ghulam Haidar\thanks{Email: \texttt{ghulam.haidar@uoh.edu.pk}}\\[1ex]
	{\small Department of Mathematics, University of Haripur, Pakistan}%
}

\date{November 2026}

\begin{document}
	\maketitle
	
	\begin{abstract}
		A connected graph is resistance-regular if the sum of the effective resistances from a vertex to all other vertices does not depend on the vertex, or equivalently if the Moore--Penrose inverse of its Laplacian has constant diagonal. Walk-regular graphs have this property, and Zhou, Wang and Bu asked in 2016 whether every resistance-regular graph is regular, a question that also arises in the study of graphs with constant resistance curvature. We answer this question in the negative and examine both what resistance-regularity forces and how far from regular such graphs can be. On the spectral side, we prove that a connected graph is walk-regular if and only if the first $m-1$ powers of its Laplacian pseudoinverse have constant diagonal, where $m$ is the number of distinct nonzero Laplacian eigenvalues. Because resistance-regularity constrains only the first power, it forces regularity when $m\le2$ and allows at most two vertex degrees when $m=3$. Any non-regular resistance-regular graph has minimum degree at least $3$ and a connected complement, and an exhaustive computer search shows that it must have at least $13$ vertices. For every even $p\ge4$, we construct a resistance-regular graph on $4p$ vertices with degrees $2p-2$ and $2p-1$, obtaining closed formulas for its Laplacian spectrum and for all pairwise effective resistances. Using module blow-ups, we then construct $n$-vertex resistance-regular graphs whose two degrees differ by $(\sqrt{4n+1}-5)/4$ for infinitely many $n$, as well as examples with three and four distinct degrees, whereas the maximum and minimum degrees of any resistance-regular graph on $n$ vertices differ by less than $n/2$. We close by conjecturing that in a resistance-regular graph the sum of squared resistances at a vertex increases strictly with its degree, which would imply that any connected graph whose first two Laplacian pseudoinverse powers both have constant diagonal is regular.
		
		
	\end{abstract}

	\section{Introduction}\label{sec:intro}
	
	When a connected graph $G$ is viewed as an electrical network in which every edge is a unit resistor, the effective resistance $\rres{u}{v}$ between two vertices $u$ and $v$ is the voltage difference induced when one unit of current enters at $u$ and leaves at $v$. Proposed as a graph metric by Klein and Randi\'c \cite{KleinRandic1993}, effective resistance connects several parts of combinatorics and probability: it is proportional to the commute time of a simple random walk \cite{Chandra1989, Tetali1991, Lovasz1993}, equals the marginal probability that an edge belongs to a uniformly random spanning tree \cite{Kirchhoff1847, DoyleSnell1984}, and underlies spectral sparsification \cite{SpielmanSrivastava2011}, current-flow centrality \cite{StephensonZelen1989, BrandesFleischer2005} and recent structural encodings for graph neural networks \cite{Zhang2023biconnectivity, Black2023oversquashing, Black2024comparing}. Expositions of the general theory can be found in Bapat's book \cite{Bapat2014} and the survey of Ellens et al.\ \cite{Ellens2011}.
	
	Write $R=(\rres{u}{v})$ for the resistance matrix of $G$, $R_v=\sum_u \rres{u}{v}$ for the resistance degree of a vertex $v$, and $\Kf(G)=\tfrac12\sum_v R_v$ for the Kirchhoff index \cite{KleinRandic1993, BapatGutmanXiao2003, Bapat2004}. Zhou, Wang and Bu \cite{ZhouWangBu2016} called a connected graph \emph{resistance-regular} if $R_v$ is independent of $v$. They showed that an $n$-vertex connected graph is resistance-regular if and only if the Moore--Penrose inverse $\Lp$ of its Laplacian $L$ has constant diagonal, and if and only if the local form of Foster's theorem \cite{Foster1949} holds at every vertex:
	\[
	\sum_{u\sim v}\rres{u}{v}=2-\frac2n\qquad\text{for all }v\in V.
	\]
	Every connected walk-regular graph in the sense of Godsil and McKay \cite{GodsilMcKay1980}, and in particular every vertex-transitive or distance-regular graph, is resistance-regular, because regularity makes $L$ a linear function of the adjacency matrix $A$, so that $\Lp$ is a polynomial in $A$ and inherits a constant diagonal. Since walk-regular graphs are necessarily regular, Zhou, Wang and Bu asked in \cite[Problem~2]{ZhouWangBu2016} whether every resistance-regular graph is regular.
	
	The same question has surfaced independently in work on discrete curvature. Devriendt and Lambiotte \cite{DevriendtLambiotte2022} defined the resistance curvature of a vertex $v$ by $p_v=1-\tfrac12\sum_{u\sim v}\rres{u}{v}$. By Foster's theorem, $\sum_v p_v=1$, so a graph has constant resistance curvature if and only if $p_v=1/n$ at every vertex, which is precisely resistance-regularity. Devriendt, Ottolini and Steinerberger \cite{DevriendtOttoliniSteinerberger2024} developed the theory of this curvature and highlighted the problem of \cite{ZhouWangBu2016}; Haritha and Chithra \cite{HarithaChithra2024} gave a characterisation of resistance-regular graphs in terms of Laplacian eigenvectors and noted that the literature contained no non-regular example; Sun, Yang and Xu \cite{SunYangXu2025} proved that graphs of constant resistance curvature are $2$-connected and that, in the normalisation $\kappa=R^{-1}\one$ of \cite{DevriendtOttoliniSteinerberger2024}, the cycle has the smallest constant curvature among them. Resistance curvature in a broader sense is studied in \cite{GuoSunYang2026}, and a related edge-based notion, in which every edge has the same effective resistance, by Palacios \cite{Palacios2026}.
	
	In this paper we show that resistance-regular graphs need not be regular, and we organise our results around three natural questions: what structural and spectral restrictions resistance-regularity imposes, how non-regular examples can be constructed, and how far their degree sequences can depart from regularity.
	
	Recall that a graph is walk-regular if the number of closed walks of each length starting at a vertex does not depend on the vertex. In terms of the local spectral measure of $L$ at a vertex $v$ (Definition~\ref{def:lsm}), the vertex degree $d_v=L_{vv}$ is the first positive moment and $\bigl((\Lp)^k\bigr)_{vv}$ is the $k$-th negative moment. Resistance-regularity requires only that the first negative moment be constant across all vertices.
	
	\begin{introtheorem}\label{thm:forced}
		Let $G$ be a connected graph whose Laplacian has $m$ distinct nonzero eigenvalues.
		\begin{enumerate}[label=(\roman*),leftmargin=2.5em]
			\item $G$ is walk-regular if and only if the diagonal entries of $\Lp,(\Lp)^2,\dots,(\Lp)^{m-1}$ are all constant. In this case the local spectral measure of $L$ is the same at every vertex, and every power of $\Lp$ has constant diagonal.
			\item If $G$ is resistance-regular and $m\le2$, then $G$ is walk-regular, and in particular regular. If $G$ is resistance-regular and $m=3$, then $G$ has at most two distinct degrees, and if it has two, they add up to one less than the sum of the three nonzero Laplacian eigenvalues.
			\item If $G$ is resistance-regular but not regular, then every vertex of $G$ has degree at least $3$; the complement of $G$ is connected, so that $G$ is not the join of two graphs and no vertex of $G$ is adjacent to all other vertices; and if $G$ is bipartite, then its two colour classes have the same size.
		\end{enumerate}
	\end{introtheorem}
	
	When $m\le2$, part (i) shows that resistance-regularity alone is strong enough to force walk-regularity, which is the mechanism underlying \cite[Theorem~14]{ZhouWangBu2016} (see Corollary~\ref{cor:three}). Once $m\ge3$, however, constancy of $\Lp_{vv}$ imposes only a single linear constraint on the $m$ weights of the local spectral measure, leaving room for the vertex degree $d_v$ to vary. Combined with the $2$-connectivity of resistance-regular graphs established by Sun, Yang and Xu \cite{SunYangXu2025}, the degree and complement restrictions in part (iii) make an exhaustive census of small graphs feasible: no non-regular resistance-regular graph exists on at most $12$ vertices, nor on $13$ vertices with maximum and minimum degrees differing by at most $2$ (Proposition~\ref{prop:enum}). Turning to graphs that admit a semiregular automorphism group with two orbits, we found sixteen non-regular examples on $16$ to $48$ vertices, the smallest of which generalises to an infinite family.
	
	\begin{introtheorem}\label{thm:existence}
		For every even integer $p\ge4$ there is a connected graph $\Fp(p)$ on $4p$ vertices, with $2p$ vertices of degree $2p-1$ and $2p$ vertices of degree $2p-2$, in which every vertex has the same resistance degree. Its Laplacian pseudoinverse has diagonal entries $(4p^2-3p-2)/(8p^2(p-2))$, its Laplacian has exactly five distinct nonzero eigenvalues, and every effective resistance in it is one of nine explicit rational functions of $p$.
	\end{introtheorem}
	
	Two of the sporadic $32$-vertex examples from the bi-Cayley search have $m=4$ (Proposition~\ref{prop:sharp}), showing that resistance-regularity fails to force regularity already for $m=4$ and $m=5$, while the case $m=3$ remains open. A key structural feature of $\Fp(p)$ is that its cliques $Q_0$ and $Q_1$ are modules, meaning that all vertices in a clique share the exact same neighbours outside it. More generally, when the vertices of a base graph are blown up into cliques, independent sets or arbitrary walk-regular modules, resistance-regularity reduces to a small system of rational equations in the module sizes. Analysing these blow-ups, alongside general Dirichlet and triangle-inequality bounds, yields our third main result.
	
	\begin{introtheorem}\label{thm:far}
		Let $\Delta$ and $\delta$ denote the maximum and minimum degrees of a connected graph $G$.
		\begin{enumerate}[label=(\roman*),leftmargin=2.5em]
			\item For every integer $t\ge2$ there is a connected resistance-regular graph on $n=4t^2+2t$ vertices with degrees $2t^2$ and $2t^2+t-1$, so that $\Delta-\delta=(\sqrt{4n+1}-5)/4$.
			\item There are connected resistance-regular graphs with exactly three and with exactly four distinct degrees.
			\item Every connected resistance-regular graph on $n$ vertices satisfies $\Delta-\delta<n/2$, and $\Delta-\delta<(2-\sqrt3)\,n$ if it has exactly two distinct degrees.
		\end{enumerate}
	\end{introtheorem}
	
	Thus the largest difference $g(n)$ between the maximum and the minimum degree of a resistance-regular graph on $n$ vertices satisfies $(\sqrt{4n+1}-5)/4\le g(n)<n/2$ for infinitely many $n$.
	
	While Theorems~\ref{thm:existence} and~\ref{thm:far} show that the row sums of the resistance matrix $R$ do not determine the vertex degrees, requiring the entire rows of $R$ to coincide as multisets is far more restrictive. In Section~\ref{sec:rows} we conjecture that any connected graph whose resistance rows are permutations of one another must be regular (Conjecture~\ref{conj:rows}), and we prove that its maximum and minimum degrees can differ by at most one (Proposition~\ref{prop:nearregular}). We also formulate a stronger monotonicity conjecture: in any resistance-regular graph, the second resistance moment $\sum_u\rres{v}{u}^2$ is a strictly increasing function of the degree $d_v$ (Conjecture~\ref{conj:moment}). We prove this monotonicity in several cases (Proposition~\ref{prop:momentcases}), verify it across all our examples, and show by a five-vertex counterexample that it fails once edges are allowed to carry unequal weights (Example~\ref{ex:weighted}). If true in general, Conjecture~\ref{conj:moment} would imply that constancy of the diagonals of both $\Lp$ and $(\Lp)^2$ suffices to force regularity for any $m$.
	
	The rest of the paper is organised as follows. Section~\ref{sec:prelim} collects the required background on $\Lp$, effective resistance and the Dirichlet lower bounds. Section~\ref{sec:forced} proves Theorem~\ref{thm:forced} and reports the exhaustive census of small graphs up to order $13$. Section~\ref{sec:construction} describes the bi-Cayley search, introduces the family $\Fp(p)$ and proves Theorem~\ref{thm:existence}. Section~\ref{sec:modules} develops the general theory of modules and blow-ups, which is then applied in Section~\ref{sec:far} to prove Theorem~\ref{thm:far} and relate the degree spread to algebraic connectivity. Section~\ref{sec:rows} studies the rows and second moments of $R$, and Section~\ref{sec:consequences} discusses implications for resistance curvature, network centrality, random spanning trees and distance-based Weisfeiler--Leman tests, concluding with open problems. Appendix~\ref{app:fourparam} analyses a four-parameter generalisation of $\Fp(p)$, and Appendix~\ref{app:examples} lists the explicit connection sets and ancillary verification files.
	
	\section{Preliminaries}\label{sec:prelim}
	
	Throughout the paper, $G=(V,E)$ denotes a finite, simple, undirected, connected graph on $n=|V|\ge3$ vertices. We write $u\sim v$ when $u$ and $v$ are adjacent, $d_v$ for the degree of $v$, $A$ for the adjacency matrix, $D=\operatorname{diag}(d_v)$ for the degree matrix, $L=D-A$ for the Laplacian, $J$ for the $n\times n$ all-ones matrix and $\one$ for the all-ones column vector.
	
	Because $G$ is connected, $L$ is positive semidefinite with a one-dimensional kernel spanned by $\one$. Its Moore--Penrose inverse $\Lp$ is the unique symmetric matrix satisfying
	\[
	L\Lp=\Lp L=I-\frac1n J\qquad\text{and}\qquad \Lp\one=0,
	\]
	and can be computed as $\Lp=(L+J/n)^{-1}-J/n$ \cite{Bapat2014}. Let $\mu_1<\mu_2<\dots<\mu_m$ be the distinct nonzero eigenvalues of $L$, and let $E_1,\dots,E_m$ be the orthogonal projections onto the corresponding eigenspaces. Then
	\[
	L=\sum_{i=1}^m \mu_i E_i,\qquad \Lp=\sum_{i=1}^m \mu_i^{-1}E_i,\qquad \sum_{i=1}^m E_i=I-\frac1n J.
	\]
	
	\begin{definition}[local spectral measure]\label{def:lsm}
		For a vertex $v\in V$, the nonnegative numbers $e_i(v)=(E_i)_{vv}$, $i=1,\dots,m$, form the \emph{local spectral measure} of $L$ at $v$. They satisfy $\sum_{i=1}^m e_i(v)=1-1/n$, and $(L^k)_{vv}=\sum_{i=1}^m \mu_i^k e_i(v)$ and $\bigl((\Lp)^k\bigr)_{vv}=\sum_{i=1}^m \mu_i^{-k} e_i(v)$ for every integer $k\ge1$.
	\end{definition}
	
	\begin{definition}[effective resistance, resistance degree, Kirchhoff index]\label{def:resistance}
		The \emph{effective resistance} between vertices $u$ and $v$ is
		\begin{equation}\label{eq:resistance}
			\rres{u}{v}=\Lp_{uu}+\Lp_{vv}-2\Lp_{uv},
		\end{equation}
		the energy dissipated by a unit current flowing from $u$ to $v$ when every edge is a unit resistor. The matrix $R=(\rres{u}{v})$ is the \emph{resistance matrix}, the \emph{resistance degree} of $v$ is $R_v=\sum_{u\in V}\rres{u}{v}$, and the \emph{Kirchhoff index} is $\Kf(G)=\tfrac12\sum_{v\in V}R_v$.
	\end{definition}
	
	Effective resistance is a metric on $V$ \cite{KleinRandic1993, Bapat2014}. Summing \eqref{eq:resistance} over $u\in V$ and using $\Lp\one=0$ relates $R_v$ to the diagonal entry $\Lp_{vv}$:
	\begin{equation}\label{eq:rowsum}
		R_v=n\Lp_{vv}+\tr\Lp.
	\end{equation}
	In particular $\Kf(G)=n\tr\Lp$. By Kirchhoff's theorem \cite{Kirchhoff1847, DoyleSnell1984}, the effective resistance across an edge $uv\in E$ equals
	\begin{equation}\label{eq:kirchhoff}
		\rres{u}{v}=\frac{\tau_{uv}}{\tau}\qquad(uv\in E),
	\end{equation}
	where $\tau$ is the number of spanning trees of $G$ and $\tau_{uv}$ is the number of spanning trees containing the edge $uv$. Summing \eqref{eq:kirchhoff} over all edges yields Foster's identity $\sum_{uv\in E}\rres{u}{v}=n-1$ \cite{Foster1949}, reflecting the fact that every spanning tree has $n-1$ edges.
	
	\begin{definition}[resistance-regular graph]\label{def:rr}
		A connected graph $G$ is \emph{resistance-regular} if the resistance degree $R_v$ is the same for all $v\in V$.
	\end{definition}
	
	We collect the basic characterisations of resistance-regularity from \cite{ZhouWangBu2016}.
	
	\begin{theorem}[{\cite[Theorem~12]{ZhouWangBu2016}}]\label{thm:zwb}
		For a connected graph $G$ on $n$ vertices, the following are equivalent:
		\begin{enumerate}[label=(\roman*),leftmargin=2.5em]
			\item $G$ is resistance-regular;
			\item the diagonal entries of $\Lp$ are all equal;
			\item $\sum_{u\sim v}\rres{u}{v}=2-2/n$ for every vertex $v\in V$;
			\item the eigenvalues of the resistance matrix $R$ are $2\Kf(G)/n$ and $-2/\mu$, where $\mu$ runs over the nonzero eigenvalues of $L$ with multiplicity.
		\end{enumerate}
	\end{theorem}
	
	The equivalence of (i) and (ii) is immediate from \eqref{eq:rowsum}, while condition (iii) states via \eqref{eq:kirchhoff} that every vertex has the same expected degree $2(n-1)/n$ in a uniformly random spanning tree. Condition (iii) also implies two simple structural restrictions that are useful when searching for examples.
	
	\begin{proposition}\label{prop:necessary}
		Let $G$ be a resistance-regular graph on $n\ge3$ vertices with $\tau$ spanning trees. Then $G$ is $2$-connected, and $n$ divides $2\tau$.
	\end{proposition}
	
	\begin{proof}
		A vertex $v$ of degree one would be incident to a single bridge $uv$, which lies in every spanning tree and therefore has $\rres{u}{v}=1$; this contradicts Theorem~\ref{thm:zwb}(iii) since $2-2/n>1$ for $n\ge3$. That $G$ has no cut vertex at all is Theorem~2.5 of \cite{SunYangXu2025}, proved there in the language of constant resistance curvature. For the divisibility condition, multiplying Theorem~\ref{thm:zwb}(iii) by $\tau$ and applying \eqref{eq:kirchhoff} gives $\sum_{u\sim v}\tau_{uv}=2\tau(n-1)/n$. The left-hand side is an integer, and since $\gcd(n,n-1)=1$, $n$ must divide $2\tau$.
	\end{proof}
	
	\begin{definition}[walk-regular graph]\label{def:walkregular}
		A graph $G$ is \emph{walk-regular} \cite{GodsilMcKay1980} if, for every integer $k\ge0$, the number of closed walks of length $k$ starting at $v$, namely $(A^k)_{vv}$, does not depend on $v$.
	\end{definition}
	
	Taking $k=2$ shows that every walk-regular graph is regular. If $G$ is connected and walk-regular of degree $d$, then $L=dI-A$ and $\Lp$ is a polynomial in $A$, so $\Lp_{vv}$ is independent of $v$ and $G$ is resistance-regular \cite[Remark~2]{ZhouWangBu2016}.
	
	\paragraph{Bounds on effective resistance.} We shall repeatedly use the following elementary lower bounds, which depend only on the degrees of the two vertices.
	
	\begin{lemma}[Dirichlet bounds]\label{lem:dirichlet}
		For any two distinct vertices $u,v$ of a connected graph $G$,
		\[
		\rres{u}{v}\ge\frac{1}{d_u}+\frac{1}{d_v}\quad\text{if }u\not\sim v,\qquad
		\rres{u}{v}\ge\frac{d_u+d_v-2}{d_ud_v-1}\quad\text{if }u\sim v,
		\]
		and in both cases $\rres{u}{v}\ge\dfrac{1}{d_u+1}+\dfrac{1}{d_v+1}$.
	\end{lemma}
	
	\begin{proof}
		By Thomson's and Dirichlet's principles, $1/\rres{u}{v}$ is the minimum of the Dirichlet energy $\sum_{ab\in E}(\psi_a-\psi_b)^2$ over all functions $\psi\colon V\to\R$ with $\psi_u=1$ and $\psi_v=0$. Restricting to test functions that take a constant value $c$ on $V\setminus\{u,v\}$ gives an upper bound on $1/\rres{u}{v}$: when $u\not\sim v$ the energy is $d_u(1-c)^2+d_vc^2$, which has minimum $d_ud_v/(d_u+d_v)$ over $c\in\R$, and when $u\sim v$ the energy is $1+(d_u-1)(1-c)^2+(d_v-1)c^2$, which has minimum $(d_ud_v-1)/(d_u+d_v-2)$. Taking reciprocals yields the first two inequalities. The unified bound follows because $\frac1{d_u}+\frac1{d_v}>\frac1{d_u+1}+\frac1{d_v+1}$ and
		\[
		\frac{d_u+d_v-2}{d_ud_v-1}-\Bigl(\frac{1}{d_u+1}+\frac{1}{d_v+1}\Bigr)=\frac{(d_u-d_v)^2}{(d_ud_v-1)(d_u+1)(d_v+1)}\ge0.\qedhere
		\]
	\end{proof}

	\section{What resistance-regularity forces}\label{sec:forced}
	
	Before constructing non-regular examples, we establish the spectral and structural constraints stated in Theorem~\ref{thm:forced} and use them to rule out small counterexamples by exhaustive enumeration.
	
	\subsection{Proof of Theorem~\ref{thm:forced}}\label{sec:forcedproof}
	
	Parts (i) and (ii) of Theorem~\ref{thm:forced} concern the moments of the local spectral measure. Write $P=I-J/n$, and for each integer $s$ define $L^{(s)}=L^s$ when $s\ge1$, $L^{(0)}=P$, and $L^{(s)}=(\Lp)^{-s}$ when $s\le-1$, so that $L^{(s)}=\sum_{i=1}^m\mu_i^sE_i$ and $\bigl(L^{(s)}\bigr)_{vv}=\sum_{i=1}^m\mu_i^se_i(v)$ for all $s\in\Z$. Alongside a Vandermonde inversion argument, the proof relies on linear relations between positive and negative powers of $L$ inherited from its minimal polynomial on $\one^\perp$.
	
	\begin{lemma}\label{lem:minpoly}
		Let $G$ be a connected graph whose Laplacian has the distinct nonzero eigenvalues $\mu_1,\dots,\mu_m$, and let $\sigma_0=1,\sigma_1,\dots,\sigma_m$ be their elementary symmetric functions. Then for every integer $j$,
		\[
		\sum_{k=0}^m(-1)^k\sigma_k\,L^{(m-k-j)}=0.
		\]
	\end{lemma}
	
	\begin{proof}
		By the spectral decomposition, the left-hand side equals $\sum_{i=1}^m c_iE_i$ with
		\[
		c_i=\sum_{k=0}^m(-1)^k\sigma_k\mu_i^{m-k-j}=\mu_i^{-j}\prod_{l=1}^m(\mu_i-\mu_l)=0.\qedhere
		\]
	\end{proof}
	
	For $m=3$ and $j=1,2$, Lemma~\ref{lem:minpoly} gives
	\begin{equation}\label{eq:m3}
		L^2-\sigma_1L+\sigma_2P-\sigma_3\Lp=0,\qquad L-\sigma_1P+\sigma_2\Lp-\sigma_3(\Lp)^2=0.
	\end{equation}
	
	\begin{proof}[Proof of Theorem~\ref{thm:forced}]
		(i) Fix a vertex $v$. Since $\sum_{i=1}^m e_i(v)=1-1/n$, the weights $e_i(v)$ of the local spectral measure at $v$ satisfy the linear system
		\begin{equation}\label{eq:vandermonde}
			\sum_{i=1}^m\mu_i^{-s}\,e_i(v)=\bigl(L^{(-s)}\bigr)_{vv},\qquad s=0,1,\dots,m-1,
		\end{equation}
		as well as the system with $\mu_i^{s}$ and $\bigl(L^{(s)}\bigr)_{vv}$ in place of $\mu_i^{-s}$ and $\bigl(L^{(-s)}\bigr)_{vv}$. The coefficient matrices of the two systems are the Vandermonde matrices in the distinct numbers $\mu_1^{-1},\dots,\mu_m^{-1}$ and $\mu_1,\dots,\mu_m$; they are invertible and do not depend on $v$. Hence the local spectral measure at $v$ is determined by the diagonal entries of $\Lp,(\Lp)^2,\dots,(\Lp)^{m-1}$ at $v$, and also by the diagonal entries of $L,L^2,\dots,L^{m-1}$ at $v$.
		
		Suppose that $G$ is walk-regular. Then $G$ is regular of some degree $d$, so $L=dI-A$, and $(L^k)_{vv}$ is a linear combination of the numbers $(A^j)_{vv}$, $0\le j\le k$, with coefficients that depend only on $d$ and $k$. Thus no diagonal entry of a power of $L$ depends on $v$, so the local spectral measure does not depend on $v$, and neither does $\bigl((\Lp)^k\bigr)_{vv}=\sum_{i=1}^m\mu_i^{-k}e_i(v)$ for any $k\ge1$. Conversely, suppose that the diagonal entries of $\Lp,(\Lp)^2,\dots,(\Lp)^{m-1}$ are all constant. Then the local spectral measure does not depend on $v$, and hence neither does $(L^k)_{vv}=\sum_{i=1}^m\mu_i^ke_i(v)$ for any $k\ge1$. Taking $k=1$ shows that $G$ is regular, say of degree $d$. Then $A=dI-L$, so $(A^k)_{vv}$ is a linear combination of the numbers $(L^j)_{vv}$, $0\le j\le k$, with coefficients that do not depend on $v$, and $G$ is walk-regular.
		
		(ii) If $m\le2$, the conditions in (i) consist of the constancy of $\Lp_{vv}$ alone, or are empty when $m=1$, so a resistance-regular graph with $m\le2$ is walk-regular by (i). Now let $m=3$. Since $(L^2)_{vv}=d_v^2+d_v$ and $P_{vv}=1-1/n$, the diagonal of the first identity in \eqref{eq:m3} reads
		\begin{equation}\label{eq:m3diag}
			\sigma_3\Lp_{vv}=d_v^2+(1-\sigma_1)d_v+\sigma_2\Bigl(1-\frac1n\Bigr).
		\end{equation}
		If $\Lp_{vv}$ is constant, every degree is therefore a root of one and the same quadratic polynomial $z^2+(1-\sigma_1)z+c$. Hence there are at most two distinct degrees, and two distinct roots add up to $\sigma_1-1$.
		
		(iii) Let $G$ be resistance-regular but not regular. By Proposition~\ref{prop:necessary}, $G$ is $2$-connected, so every vertex has degree at least $2$. Suppose for a contradiction that some vertex $v$ has degree $2$, with neighbours $a$ and $b$. Since $G$ is $2$-connected, $G'=G-v$ is connected on $n-1$ vertices; let $r'$ denote the effective resistance between $a$ and $b$ in $G'$. Because $va$ and $vb$ are the only edges incident to $v$, the network between $v$ and $a$ places the unit resistor $va$ in parallel with the series combination of $vb$ and $G'$. Thus $\rres{v}{a}=(1+r')/(2+r')$, and by symmetry $\rres{v}{b}=(1+r')/(2+r')$ as well. Applying Theorem~\ref{thm:zwb}(iii) at $v$ gives $2(1+r')/(2+r')=2-2/n$, which forces $r'=n-2$. By Thomson's principle \cite{DoyleSnell1984}, however, $r'$ is the minimum energy of a unit flow from $a$ to $b$ in $G'$ and therefore cannot exceed the length of a shortest $a$--$b$ path $P'$ in $G'$, which is at most $n-2$. Equality requires $P'$ to have length $n-2$ (visiting all $n-1$ vertices of $G'$) and the unit flow along $P'$ to be the unique electrical current from $a$ to $b$ in $G'$. Since the potential drops by $1$ along each edge of $P'$, any additional edge in $G'$ would connect two vertices at different potentials and carry nonzero current. It follows that $G'=P'$ and $G$ is the cycle $C_n$, contradicting the assumption that $G$ is non-regular. Every vertex of $G$ must therefore have degree at least $3$.
		
		Next, let $\overline G$ be the complement of $G$, with Laplacian $\overline L$ and degrees $h_v=n-1-d_v$. Since $L=nI-J-\overline L$ and $\overline L\one=0$, the matrix $nI-\overline L$ is invertible (its eigenvalue on $\one$ is $n$, and on $\one^\perp$ it agrees with $L$) and $\Lp=(nI-\overline L)^{-1}-J/n^2$. Hence $G$ is resistance-regular if and only if $f_v=\bigl[(nI-\overline L)^{-1}\bigr]_{vv}$ is independent of $v$. Let $\nu_v$ be the local spectral measure of $\overline L$ at $v$, now including the eigenvalue $0$, so that $\nu_v$ is a probability measure with mean $\int\lambda\,d\nu_v=h_v$ and second moment $(\overline L^2)_{vv}=h_v^2+h_v$, hence variance $h_v$. If $v$ lies in a component of $\overline G$ with $c_v$ vertices, $\nu_v$ is supported on the Laplacian eigenvalues of that component, which lie in $[0,c_v]$ \cite{Mohar1991}. Integrating the identity
		\[
		\frac1{n-\lambda}=\frac1{n-h}+\frac{\lambda-h}{(n-h)^2}+\frac{(\lambda-h)^2}{(n-h)^2(n-\lambda)}
		\]
		with $h=h_v$ against $\nu_v$ gives $f_v=\frac1{n-h_v}+\frac{K_v}{(n-h_v)^2}$ with $0\le K_v\le h_v/(n-c_v)$ whenever $c_v<n$. Now suppose that $\overline G$ is disconnected and that $u$ and $w$ lie in different components with $d_u<d_w$, that is, $h_u>h_w$. Then $h_w<h_u\le c_u-1\le n-c_w-1$, so $K_w<1$, and
		\[
		f_u-f_w\ge\frac1{n-h_u}-\frac1{n-h_w}-\frac{K_w}{(n-h_w)^2}=\frac{(h_u-h_w)(n-h_w)-K_w(n-h_u)}{(n-h_u)(n-h_w)^2}>0,
		\]
		because $h_u-h_w\ge1$ and $n-h_w>n-h_u$, which contradicts resistance-regularity. Thus, if $\overline G$ is disconnected, any two vertices in different components of $\overline G$ have the same degree in $G$; since every vertex has a vertex in another component, all degrees are equal, contrary to our assumption. Hence $\overline G$ is connected. In particular, $G$ is not the join of two graphs, and no vertex of $G$ is adjacent to all other vertices, since such a vertex would be isolated in $\overline G$.
		
		Finally, every edge of a bipartite graph has exactly one end in each colour class $U$. If $G$ is bipartite, summing Theorem~\ref{thm:zwb}(iii) over $v\in U$ and using Foster's identity gives $|U|(2-2/n)=\sum_{uv\in E}\rres{u}{v}=n-1$, so $|U|=n/2$.
	\end{proof}
	
	When $m\ge3$, the two equations $s=0$ and $s=1$ of \eqref{eq:vandermonde} no longer determine the $m$ weights $e_i(v)$, and Section~\ref{sec:construction} shows that the degrees can then indeed differ. For $m=2$, by contrast, Lemma~\ref{lem:minpoly} with $j=1$ gives $L-\sigma_1P+\sigma_2\Lp=0$, whose diagonal $d_v=\sigma_1(1-1/n)-\sigma_2\Lp_{vv}$ expresses the degree through $\Lp_{vv}$. Zhou, Wang and Bu \cite[Theorem~14]{ZhouWangBu2016} derived this identity for a graph with two distinct nonzero Laplacian eigenvalues and then concluded that $G$ is regular from the constancy of $\Lp_{vv}$; reading the diagonal of the identity is what makes the conclusion valid. Since constancy of $\Lp_{vv}$ alone does not imply regularity (Theorem~\ref{thm:existence}), a linear relation of this kind between $L$ and $\Lp$, which exists exactly when $m=2$, is essential to their argument. We record the resulting statement.
	
	\begin{corollary}[{\cite[Theorem~14]{ZhouWangBu2016}}]\label{cor:three}
		A connected resistance-regular graph whose resistance matrix $R$ has exactly three distinct eigenvalues is strongly regular.
	\end{corollary}
	
	\begin{proof}
		By Theorem~\ref{thm:zwb}(iv), the eigenvalues of $R$ are $2\Kf(G)/n>0$ together with $-2/\mu_i<0$ for $1\le i\le m$. Thus $R$ has exactly three distinct eigenvalues if and only if $L$ has $m=2$ distinct nonzero eigenvalues. By Theorem~\ref{thm:forced}(ii), $G$ is regular, and a connected regular graph with three distinct eigenvalues is strongly regular \cite[Lemma~10.2.1]{GodsilRoyle2001}.
	\end{proof}
	
	The second identity in \eqref{eq:m3} also controls the next negative moment when $m=3$; we shall use this in Section~\ref{sec:moment}.
	
	\begin{corollary}\label{cor:m3}
		Let $G$ be a connected resistance-regular graph whose Laplacian has exactly three distinct nonzero eigenvalues. Then $\bigl((\Lp)^2\bigr)_{vv}$ is a strictly increasing function of $d_v$.
	\end{corollary}
	
	\begin{proof}
		The diagonal of the second identity in \eqref{eq:m3} reads $\sigma_3\bigl((\Lp)^2\bigr)_{vv}=d_v-\sigma_1(1-1/n)+\sigma_2\Lp_{vv}$. Since $\Lp_{vv}$ is constant and $\sigma_3>0$, the claim follows.
	\end{proof}
	
	\subsection{Small graphs}\label{sec:enum}
	
	Theorem~\ref{thm:forced}(iii) and Proposition~\ref{prop:necessary} make an exhaustive search of small graphs feasible: a non-regular resistance-regular graph is $2$-connected and has minimum degree at least $3$. We enumerated graphs with \texttt{geng} from the \texttt{nauty} package \cite{McKayPiperno2014} and used two tests (Table~\ref{tab:enum}). For $n\le10$ we tested all connected graphs in exact integer arithmetic: fraction-free Gaussian elimination of the Laplacian with one row and column removed gives the number $\tau$ of spanning trees and the integers $\tau\rres{u}{v}$, and a graph is resistance-regular exactly when the integer row sums $\sum_{u}\tau\rres{u}{v}$ are equal.
	
	For $11\le n\le13$ we restricted the enumeration to $2$-connected graphs with minimum degree at least $3$ and compiled a two-stage floating-point filter directly into \texttt{geng} so that graphs could be screened on the fly without disk storage. In the first stage, the diagonal entries of $(L+J/n)^{-1}$ at one vertex of maximum degree and one of minimum degree are estimated in single precision by eliminating the other $n-2$ vertices, discarding the graph only if those two entries differ by more than $10\%$. For every graph that survives this coarse check, the second stage computes the full diagonal of $(L+J/n)^{-1}$ in double precision via Cholesky factorisation and flags the graph as a candidate whenever the relative spread $(\max-\min)/\min$ falls below $10^{-7}$. Both stages rely solely on pivot-free symmetric elimination of the positive definite matrix $L+J/n$, whose eigenvalues lie in $[\lambda,n]$ with $\lambda\ge\min\{1,4(1-\cos(\pi/n))\}$ for $2$-connected graphs \cite{Fiedler1973}. Standard backward-error analysis \cite{Higham2002} guarantees that for $n\le13$ the relative error in every computed diagonal entry is below $4\times10^{-3}$ in the first stage and below $10^{-10}$ in the second (see the ancillary documentation for details), leaving a wide safety margin below both thresholds. As an empirical validation, the compiled filter correctly identified random relabellings of all $23$ known non-regular resistance-regular graphs on at most $32$ vertices, and in the regular degree classes tested on $12$ and $13$ vertices (degrees $3,4,5$ for $n=12$ and degrees $4,6$ for $n=13$) it accepted the exact set of resistance-regular graphs determined by rational arithmetic.
	
	For $n=12$ this covers all $113\,372\,069\,299$ two-connected graphs of minimum degree at least $3$, of which $18\,972$ are regular. No candidate appeared; the smallest relative spread among the non-regular graphs is $4/125$, attained by a graph with six vertices of degree $4$ and six of degree $5$. For $n=13$, where there are about $5\times10^{13}$ connected graphs, we covered all $2$-connected graphs of minimum degree at least $3$ whose degrees lie in a window $\{k,k+1,k+2\}$ for $k=3,\dots,7$, and all graphs of minimum degree at least $8$. These runs include every graph on $13$ vertices whose maximum and minimum degree differ by at most $2$ and that is not excluded by Proposition~\ref{prop:necessary} and Theorem~\ref{thm:forced}(iii). The runs were split into residue classes of \texttt{geng}, and the coverage and the totals in Table~\ref{tab:enum} can be checked from the per-chunk records in the ancillary files. As an independent check, a double-precision screen of all connected graphs of minimum degree at least $2$ on at most $10$ vertices and of all $900\,969\,091$ two-connected graphs on $11$ vertices (a number that agrees with sequence A002218 in \cite{OEIS}), with exact verification of every graph passing the screen, found no non-regular resistance-regular graph either.
	
	\begin{table}[!htbp]
		\centering
		\caption{Exhaustive enumerations of small graphs with \texttt{geng}; $\delta$ denotes the minimum degree; \emph{exact} refers to the test in exact integer arithmetic and \emph{built-in} to the floating-point test built into \texttt{geng}, both described in the text. Every resistance-regular graph in these classes is regular. For $n=13$ the counts are summed over the degree windows, so graphs whose degrees lie in two windows are counted twice.}
		\label{tab:enum}
		\begin{adjustbox}{max width=\textwidth}
			\begin{tabular}{@{}llllr@{}}
				\toprule
				order $n$ & graph class enumerated & degree restriction & test & graphs tested\\
				\midrule
				$3$--$10$ & all connected graphs & none & exact & $11\,989\,762$\\
				$11$ & $2$-connected, $\delta\ge3$ & none & built-in & $577\,050\,233$\\
				$12$ & $2$-connected, $\delta\ge3$ & none & built-in & $113\,372\,069\,299$\\
				$13$ & $2$-connected, $\delta\ge3$ & degrees in $\{k,k+1,k+2\}$, $k=3,\dots,7$ & built-in & $152\,893\,105\,672$\\
				$13$ & all graphs & $\delta\ge8$ & built-in & $77\,423\,446$\\
				\bottomrule
			\end{tabular}
		\end{adjustbox}
	\end{table}
	
	\begin{proposition}[computer-assisted]\label{prop:enum}
		Every resistance-regular graph on at most $12$ vertices is regular. Every resistance-regular graph on $13$ vertices whose maximum and minimum degree differ by at most $2$, or whose minimum degree is at least $8$, is regular.
	\end{proposition}
	
	\begin{proof}
		No non-regular resistance-regular graph was found in any of the enumerations summarised in Table~\ref{tab:enum}. By Proposition~\ref{prop:necessary} and Theorem~\ref{thm:forced}(iii), any graph on $11$ or $12$ vertices omitted from the search, as well as any $13$-vertex graph satisfying the stated degree conditions but omitted from the search, cannot be a non-regular resistance-regular graph.
	\end{proof}
	
	Consequently, the minimum order of a non-regular resistance-regular graph is at least $13$; the graph $\Fp(4)$ constructed in the next section shows that it is at most $16$.

	\section{Existence: a family of non-regular resistance-regular graphs}\label{sec:construction}
	
	By Proposition~\ref{prop:enum}, a non-regular resistance-regular graph has at least $13$ vertices, and by Theorem~\ref{thm:forced}(ii) its Laplacian has at least three distinct nonzero eigenvalues. Since an exhaustive search becomes infeasible beyond $13$ vertices, we turned to graphs with a large group of automorphisms, for which the diagonal of $\Lp$ can be computed from small matrices. In this section we describe this search and the infinite family that it revealed, and we prove Theorem~\ref{thm:existence}.
	
	\subsection{Bi-Cayley graphs and a computer search}\label{sec:search}
	
	We consider graphs that have a semiregular group of automorphisms with two orbits; these are the bi-Cayley graphs, whose automorphisms were studied by Zhou and Feng \cite{ZhouFeng2016}.
	
	\begin{definition}[bi-Cayley graph]\label{def:bicayley}
		Let $H$ be a finite abelian group, written additively, with character group $\widehat H$. Given inverse-closed subsets $S_{11},S_{22}\subseteq H\setminus\{0\}$ and an arbitrary subset $S_{12}\subseteq H$, the \emph{bi-Cayley graph} $\mathrm{BiCay}(H;S_{11},S_{22},S_{12})$ has vertex set $H\times\{1,2\}$ and edges
		\[
		(h,1)\sim(h+s,1)\ (s\in S_{11}),\qquad (h,2)\sim(h+s,2)\ (s\in S_{22}),\qquad (h,1)\sim(h+s,2)\ (s\in S_{12}).
		\]
	\end{definition}
	
	Left translation by $H$ acts semiregularly on the vertices with the two orbits $H\times\{1\}$ and $H\times\{2\}$, which have constant degrees $d_1=|S_{11}|+|S_{12}|$ and $d_2=|S_{22}|+|S_{12}|$. For a character $\chi\in\widehat H$ and a subset $S\subseteq H$, write $\chi(S)=\sum_{s\in S}\chi(s)$; when $S=-S$, $\chi(S)$ is real.
	
	Because the four orbit blocks of $L$ lie in the group algebra of $H$, the characters of $H$ simultaneously block-diagonalise $L$. Specifically, $L$ is unitarily similar to the direct sum over $\chi\in\widehat H$ of the $2\times2$ Hermitian matrices
	\begin{equation}\label{eq:block}
		\widehat L_\chi=\begin{pmatrix} d_1-\chi(S_{11}) & -\chi(S_{12})\\ -\overline{\chi(S_{12})} & d_2-\chi(S_{22})\end{pmatrix}
	\end{equation}
	(see \cite{GaoLuo2010}). Inverting these $2\times2$ blocks gives an explicit criterion for resistance-regularity.
	
	\begin{lemma}\label{lem:bicayley}
		Let $G=\mathrm{BiCay}(H;S_{11},S_{22},S_{12})$ with $H$ abelian and $S_{12}\ne\emptyset$. Then $G$ is connected if and only if $\det\widehat L_\chi\ne0$ for every nontrivial character $\chi\in\widehat H\setminus\{1\}$. In that case, for every $h\in H$,
		\begin{align}
			\Lp_{(h,1),(h,1)}&=\frac{1}{4|H|\,|S_{12}|}+\frac{1}{|H|}\sum_{\chi\ne1}\frac{d_2-\chi(S_{22})}{\det\widehat L_\chi},\label{eq:a1}\\
			\Lp_{(h,2),(h,2)}&=\frac{1}{4|H|\,|S_{12}|}+\frac{1}{|H|}\sum_{\chi\ne1}\frac{d_1-\chi(S_{11})}{\det\widehat L_\chi}.\label{eq:a2}
		\end{align}
		In particular, $G$ is resistance-regular if and only if
		\begin{equation}\label{eq:criterion}
			\sum_{\chi\ne1}\frac{\bigl(d_2-\chi(S_{22})\bigr)-\bigl(d_1-\chi(S_{11})\bigr)}{\det\widehat L_\chi}=0.
		\end{equation}
	\end{lemma}
	
	\begin{proof}
		Let $U$ be the unitary character table of $H$ with columns $|H|^{-1/2}\chi$ for $\chi\in\widehat H$. Conjugating $L$ by $U\oplus U$ and reordering coordinates block-diagonalises $L$ into $\bigoplus_{\chi\in\widehat H}\widehat L_\chi$, so $\dim\ker L=\sum_\chi\dim\ker\widehat L_\chi$. For the trivial character $\chi=1$, we have $\chi(S_{ii})=|S_{ii}|$ and $\chi(S_{12})=|S_{12}|$, so $\widehat L_1=|S_{12}|\bigl(\begin{smallmatrix}1&-1\\-1&1\end{smallmatrix}\bigr)$, which has a one-dimensional kernel because $S_{12}\ne\emptyset$. Thus $\ker L$ is one-dimensional, meaning that $G$ is connected, if and only if $\det\widehat L_\chi\ne0$ for all $\chi\ne1$.
		
		Since the Moore--Penrose inverse respects unitary similarity and direct sums, $\Lp=(U\oplus U)\bigl(\bigoplus_\chi\widehat L_\chi^{\dagger}\bigr)(U\oplus U)^*$. Because every entry of $U$ has modulus $|H|^{-1/2}$, the diagonal entry of $\Lp$ at $(h,i)$ is simply the average of the $(i,i)$ entries of $\widehat L_\chi^{\dagger}$ over all $\chi\in\widehat H$. For $\chi=1$, $\widehat L_1^{\dagger}=(4|S_{12}|)^{-1}\bigl(\begin{smallmatrix}1&-1\\-1&1\end{smallmatrix}\bigr)$, which contributes $(4|H|\,|S_{12}|)^{-1}$ to both \eqref{eq:a1} and \eqref{eq:a2}. For $\chi\ne1$, $\widehat L_\chi$ is invertible and Cramer's rule gives $(\widehat L_\chi^{-1})_{11}=(d_2-\chi(S_{22}))/\det\widehat L_\chi$ and $(\widehat L_\chi^{-1})_{22}=(d_1-\chi(S_{11}))/\det\widehat L_\chi$. This establishes \eqref{eq:a1} and \eqref{eq:a2}. Since $H$ acts transitively on each of the two orbits, $\Lp$ has constant diagonal if and only if \eqref{eq:a1} and \eqref{eq:a2} are equal, which is \eqref{eq:criterion}.
	\end{proof}

	An analogous decomposition holds for non-abelian groups $H$, with characters replaced by irreducible unitary representations $\rho$ and the $2\times2$ blocks \eqref{eq:block} replaced by blocks of size $2\dim\rho$; we used this general form in the search described next.
	
	Lemma~\ref{lem:bicayley} allows us to test a bi-Cayley graph for resistance-regularity by evaluating $2\times2$ determinants rather than inverting an $n\times n$ matrix. For each finite group $H$, we enumerated triples $(S_{11},S_{22},S_{12})$ of subsets of $H$ such that $S_{11},S_{22}\subseteq H\setminus\{e\}$ are inverse-closed with $|S_{11}|<|S_{22}|$ and $e\in S_{12}$ (both normalisations are without loss of generality up to isomorphism). Keeping the connected graphs with $d_1\ne d_2$, we screened each candidate via \eqref{eq:criterion}, or its representation-theoretic analogue for non-abelian $H$, in double precision and verified every surviving graph in exact rational arithmetic. When $H$ is abelian, the left-hand side of \eqref{eq:criterion} depends on $S_{12}$ only through the squared moduli $|\chi(S_{12})|^2$, which makes it possible to batch all subsets $S_{12}$ sharing the same Fourier profile. Using this reduction, we searched exhaustively over all groups of order at most $12$, all cyclic groups $\Z_k$ with $k\le20$, and the order-$16$ abelian groups $\Z_{16}$, $\Z_2\times\Z_8$, $\Z_4\times\Z_4$ and $\Z_2^2\times\Z_4$ ($364\,456\,403\,828$ connected labelled triples in total before Fourier grouping), and we also sampled the dihedral groups $D_7,\dots,D_{10}$. Passing graphs were canonicalised up to isomorphism using \texttt{nauty} \cite{McKayPiperno2014}.
	
	\begin{table}[!htbp]
		\centering
		\caption{The sixteen non-regular resistance-regular graphs found in the bi-Cayley search, ordered by $n$ and $|E|$; $m$ is the number of distinct nonzero Laplacian eigenvalues. Graphs~1, 4 and~13 are $\Fp(4)$, $\Fp(6)$ and $\Fp(8)$. Connection sets, and the groups over which each graph arises, are given in Appendix~\ref{app:examples}.}
		\label{tab:examples}
		\begin{adjustbox}{max width=\textwidth}
			\begin{tabular}{@{}r r r l l r @{\hspace{2.5em}} r r r l l r@{}}
				\toprule
				no. & $n$ & $|E|$ & degrees & $\Lp_{vv}$ & $m$ & no. & $n$ & $|E|$ & degrees & $\Lp_{vv}$ & $m$\\
				\midrule
				1 & 16 & 52 & $6,7$ & $25/128$ & 5 & 9 & 32 & 168 & $10,11$ & $13/128$ & 5\\
				2 & 24 & 78 & $6,7$ & $121/576$ & 5 & 10 & 32 & 168 & $10,11$ & $1325/13056$ & 6\\
				3 & 24 & 78 & $6,7$ & $2857/15840$ & 10 & 11 & 32 & 184 & $11,12$ & $2555/26112$ & 6\\
				4 & 24 & 126 & $10,11$ & $31/288$ & 5 & 12 & 32 & 216 & $13,14$ & $303/3584$ & 6\\
				5 & 32 & 120 & $7,8$ & $67/448$ & 6 & 13 & 32 & 232 & $14,15$ & $115/1536$ & 5\\
				6 & 32 & 152 & $9,10$ & $907/7680$ & 8 & 14 & 32 & 232 & $14,15$ & $37/512$ & 4\\
				7 & 32 & 168 & $10,11$ & $13/128$ & 5 & 15 & 32 & 232 & $14,15$ & $37/512$ & 4\\
				8 & 32 & 168 & $10,11$ & $13/128$ & 5 & 16 & 48 & 348 & $14,15$ & $373/5040$ & 9\\
				\bottomrule
			\end{tabular}
		\end{adjustbox}
	\end{table}
	
	The search produced the sixteen pairwise non-isomorphic graphs of Table~\ref{tab:examples}, each with two consecutive degrees carried by the two orbits of $H$. Over the groups of order at most $8$, graph~1, on $16$ vertices, is the only example. Graphs~1, 4 and~13 belong to the family of Section~\ref{sec:family}, graph~12 is obtained from graph~1 by a construction of Section~\ref{sec:modules}, and we regard the other twelve as sporadic, although some of them are blow-ups in the sense of that section. Graphs~7, 8 and~9 are cospectral but pairwise non-isomorphic, as are graphs~14 and~15. Graph~3 is a bi-circulant over $\Z_{12}$: its first orbit induces the square of the $12$-cycle, its second orbit two disjoint copies of $K_6$ on the even and the odd residues, and vertex $i$ of the first orbit is joined to the vertices $i$ and $i+6$ of the second. Graph~16 has an analogous structure over $\Z_{24}$, but the same pattern gives no further examples over $\Z_k$ for $k\le40$.
	
	\subsection{The family \texorpdfstring{$\Fp(p)$}{F(p)}}\label{sec:family}
	
	Graph~1 of Table~\ref{tab:examples} has $16$ vertices of degrees $6$ and $7$. Inspecting its structure shows that its $16$ vertices partition into two $4$-cliques $Q_0,Q_1$ of degree $7$ and two outer sets $X_0,X_1$ of degree $6$: each clique $Q_b$ is completely joined to $X_b$, while the eight vertices of $X_0\cup X_1$ partition into two induced $4$-cycles that weave between the two blocks. Viewing a $4$-cycle as the cocktail-party graph $K_4$ minus a perfect matching suggests replacing $K_4$ by $K_p$ for any even integer $p\ge4$. As we now verify, this generalisation preserves the exact cancellation between the two vertex orbits for every even $p\ge4$.
	
	In the rest of this section, $p\ge4$ is an even integer and $h=p/2$.
	
	\begin{definition}\label{def:family}
		The graph $\Fp(p)$ has vertex set $Q_0\cup Q_1\cup X_0\cup X_1$, where $Q_0$ and $Q_1$ are disjoint sets of size $p$ and $X_b=\{x_{b,c,t}: c\in\{0,1\},\,t\in\Z_h\}$ for $b\in\{0,1\}$. Its edges are of three types:
		\begin{enumerate}[label=(E\arabic*),leftmargin=3em]
			\item every pair of distinct vertices inside $Q_0$ and every pair inside $Q_1$;
			\item every pair $\{q,x\}$ with $q\in Q_b$ and $x\in X_b$, for $b\in\{0,1\}$;
			\item for each $c\in\{0,1\}$, every pair of distinct vertices in $Y_c=\{x_{b,c,t}:b\in\{0,1\},\,t\in\Z_h\}$ except the $h$ matching pairs $\{x_{0,c,t},x_{1,c,t}\}$ ($t\in\Z_h$).
		\end{enumerate}
	\end{definition}
	
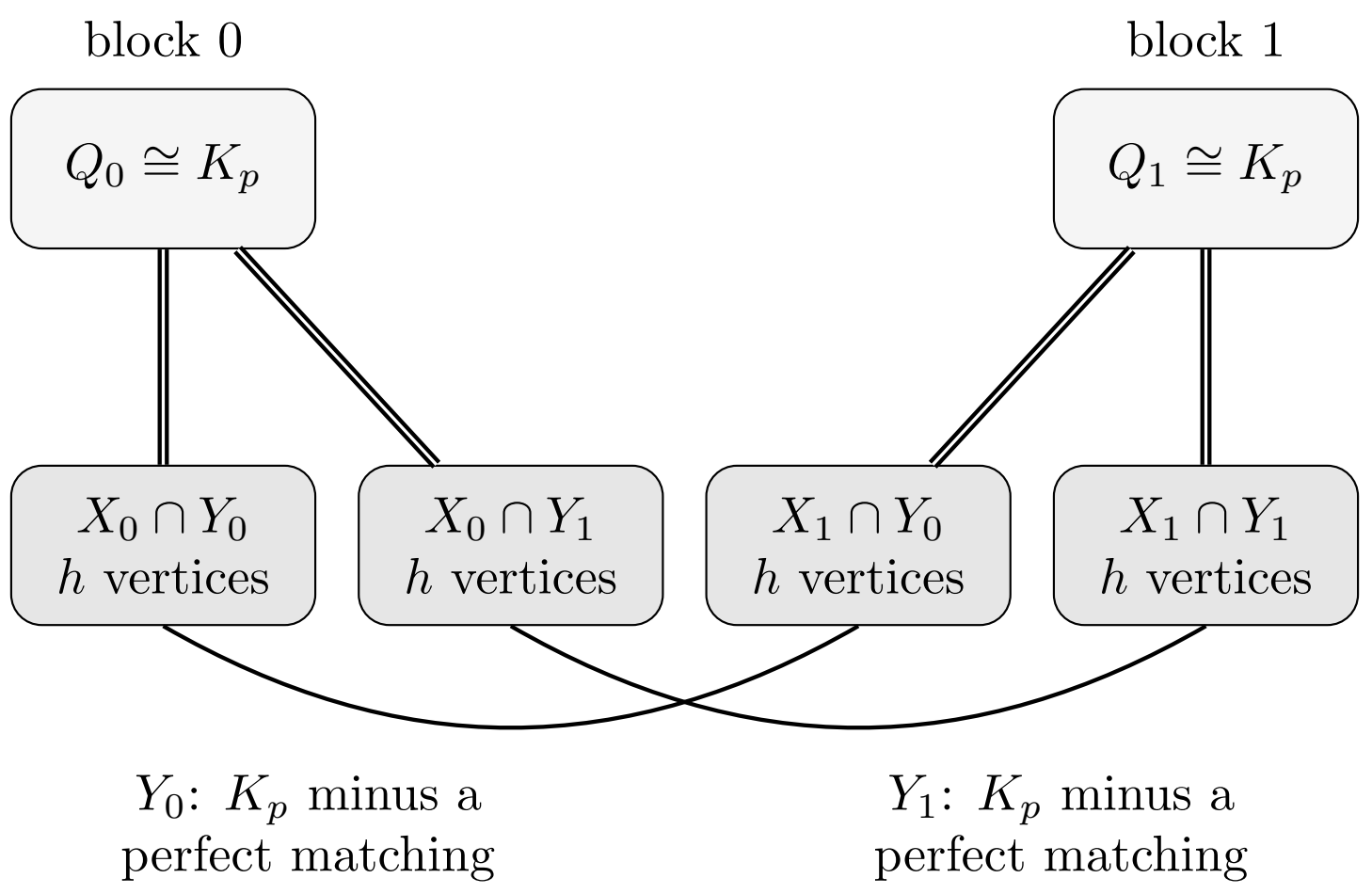
\begin{figure}[!htbp]
		\centering
		\begin{adjustbox}{width=0.6\textwidth}
			\begin{tikzpicture}[scale=1, every node/.style={font=\small}]
				\tikzset{blob/.style={draw, rounded corners=6pt, minimum width=2.1cm, minimum height=1.1cm, align=center, fill=gray!8}}
				\node[blob] (Q0) at (-3.6,2.2) {$Q_0\cong K_p$};
				\node[blob] (Q1) at (3.6,2.2) {$Q_1\cong K_p$};
				\node[blob, fill=gray!20] (X0a) at (-3.6,-0.4) {$X_0\cap Y_0$\\ $h$ vertices};
				\node[blob, fill=gray!20] (X0b) at (-1.2,-0.4) {$X_0\cap Y_1$\\ $h$ vertices};
				\node[blob, fill=gray!20] (X1a) at (1.2,-0.4) {$X_1\cap Y_0$\\ $h$ vertices};
				\node[blob, fill=gray!20] (X1b) at (3.6,-0.4) {$X_1\cap Y_1$\\ $h$ vertices};
				\draw[double, thick] (Q0) -- (X0a);
				\draw[double, thick] (Q0) -- (X0b);
				\draw[double, thick] (Q1) -- (X1a);
				\draw[double, thick] (Q1) -- (X1b);
				\draw[thick, bend right=30] (X0a.south) to (X1a.south);
				\draw[thick, bend right=30] (X0b.south) to (X1b.south);
				\node[align=center] at (-2.6,-2.35) {$Y_0$: $K_p$ minus a\\ perfect matching};
				\node[align=center] at (2.6,-2.35) {$Y_1$: $K_p$ minus a\\ perfect matching};
				\node at (-3.6,3.1) {block $0$};
				\node at (3.6,3.1) {block $1$};
			\end{tikzpicture}
		\end{adjustbox}
		\caption{The structure of $\Fp(p)$. Double lines indicate complete joins. Inside $Y_c$ every two vertices are adjacent except for the $h$ pairs that lie in different blocks and share the same index $t$.}
		\label{fig:family}
	\end{figure}

	As illustrated in Figure~\ref{fig:family}, $\Fp(p)$ consists of two symmetric blocks indexed by $b\in\{0,1\}$. Block $b$ contains an internal clique $Q_b\cong K_p$ joined completely to $p$ outer vertices $X_b$. Across the two blocks, the $2p$ outer vertices split into two sets $Y_0$ and $Y_1$ of size $p$, each taking $h$ vertices from $X_0$ and $h$ from $X_1$. Each $Y_c$ induces a cocktail-party graph ($K_p$ minus a perfect matching) in which every missing edge pairs a vertex of $X_0$ with the corresponding vertex of $X_1$. Every vertex in $Q_b$ has $p-1$ neighbours in $Q_b$ and $p$ in $X_b$, giving degree $2p-1$; every vertex in $X_b$ has $p$ neighbours in $Q_b$ and $p-2$ in its $Y_c$, giving degree $2p-2$. Thus $\Fp(p)$ is a non-regular graph on $n=4p$ vertices with $p(4p-3)$ edges.
	
	We first realise $\Fp(p)$ as a bi-Cayley graph and evaluate the character sums that enter Lemma~\ref{lem:bicayley}.
	
	\begin{lemma}\label{lem:Fbicayley}
		Let $p\ge4$ be even, $h=p/2$, and $H=\Z_2\times\Z_2\times\Z_h$, with elements written as triples $(c,b,t)$. Let $K=\Z_2\times\{0\}\times\Z_h$ be the index-$2$ subgroup of elements with $b=0$, and put
		\[
		S_{12}=K,\qquad S_{22}=K\setminus\{0\},\qquad S_{11}=\{(0,b,t)\in H: (b,t)\notin\{(0,0),(1,0)\}\}.
		\]
		Then $\Fp(p)\cong\mathrm{BiCay}(H;S_{11},S_{22},S_{12})$, with the first orbit corresponding to $X_0\cup X_1$ and the second to $Q_0\cup Q_1$, so that $d_1=2p-2$ and $d_2=2p-1$. With $M_{11}=d_1-\chi(S_{11})$, $M_{22}=d_2-\chi(S_{22})$ and $\Theta=\det\widehat L_\chi=M_{11}M_{22}-|\chi(S_{12})|^2$, the $2p-1$ nontrivial characters $\chi$ of $H$ fall into the following five classes:
		\begin{center}
			\renewcommand{\arraystretch}{1.25}
			\begin{tabular}{@{}lcccccc@{}}
				\toprule
				character class & count & $|\chi(S_{12})|$ & $M_{11}$ & $M_{22}$ & $\Theta$ & $(M_{22}-M_{11})/\Theta$\\
				\midrule
				$(\alpha,\gamma,k)=(0,1,0)$ & $1$ & $p$ & $2p-2$ & $p$ & $p(p-2)$ & $-1/p$\\
				$(1,0,0)$ & $1$ & $0$ & $p$ & $2p$ & $2p^2$ & $1/(2p)$\\
				$(1,1,0)$ & $1$ & $0$ & $2p-2$ & $2p$ & $4p(p-1)$ & $1/(2p(p-1))$\\
				$k\ne0,\ \gamma=0$ & $p-2$ & $0$ & $2p$ & $2p$ & $4p^2$ & $0$\\
				$k\ne0,\ \gamma=1$ & $p-2$ & $0$ & $2p-2$ & $2p$ & $4p(p-1)$ & $1/(2p(p-1))$\\
				\bottomrule
			\end{tabular}
		\end{center}
		Here $\chi_{\alpha,\gamma,k}(c,b,t)=(-1)^{\alpha c}(-1)^{\gamma b}\omega^{kt}$ for $\alpha,\gamma\in\Z_2$ and $k\in\Z_h$, where $\omega=e^{2\pi i/h}$.
	\end{lemma}
	
	\begin{proof}
		Both $S_{11}$ and $S_{22}$ are closed under negation and exclude $0$. Identify the first orbit $H\times\{1\}$ with $X_0\cup X_1$ by sending $((c,b,t),1)$ to $x_{b,c,t}$, and identify the second orbit $H\times\{2\}$ with $Q_0\cup Q_1$ by using $(c,t)\in\Z_2\times\Z_h$ to index the $p$ vertices of $Q_b$. Under this identification:
		\begin{itemize}[leftmargin=2em,itemsep=1pt]
			\item Two vertices in $H\times\{2\}$ are adjacent if and only if their difference lies in $S_{22}=K\setminus\{0\}$, which means they are distinct and have the same $b$ coordinate; this gives the cliques $Q_0$ and $Q_1$ in (E1).
			\item A vertex $((c,b,t),1)$ is adjacent to $((c,b,t)+k,2)$ for all $k\in S_{12}=K$, which ranges over all $p$ vertices of $Q_b$; this gives the complete joins in (E2).
			\item Two vertices in $H\times\{1\}$ are adjacent if and only if their difference lies in $S_{11}$, meaning they share the same $c$ coordinate and are neither equal nor a matching pair $\{x_{0,c,t},x_{1,c,t}\}$; this gives the cocktail-party subgraphs $Y_0$ and $Y_1$ in (E3).
		\end{itemize}
		Hence $\Fp(p)\cong\mathrm{BiCay}(H;S_{11},S_{22},S_{12})$ with orbit degrees $d_1=|S_{11}|+|S_{12}|=(p-2)+p=2p-2$ and $d_2=|S_{22}|+|S_{12}|=(p-1)+p=2p-1$. The characters of $H$ are the maps $\chi_{\alpha,\gamma,k}$, and summing $\chi=\chi_{\alpha,\gamma,k}$ over each connection set gives
		\[
		\chi(S_{12})=p\,[\alpha=0,\,k=0],\qquad
		\chi(S_{22})=\chi(S_{12})-1,\qquad
		\chi(S_{11})=p\,[\gamma=0,\,k=0]-1-(-1)^{\gamma},
		\]
		where $[\cdot]$ denotes the Iverson bracket. Evaluating $M_{11}$, $M_{22}$ and $\Theta$ for the nontrivial characters gives the table.
	\end{proof}
	
	The character sums of Lemma~\ref{lem:Fbicayley} lead quickly to $\Lp_{vv}$, as the proof of Theorem~\ref{thm:existence} below shows, but solving the potential equations directly by symmetry gives every entry of $\Lp$ and shows why the edge resistances at the two vertex orbits balance. For a fixed vertex $v$, the column $\varphi=\Lp e_v$ is the unique vector satisfying $L\varphi=e_v-\one/n$ and $\one^{\top}\varphi=0$; physically, $\varphi$ is the vertex potential when one unit of current enters at $v$ and $1/n$ leaves at every vertex. Because any automorphism of $\Fp(p)$ fixing $v$ preserves $\varphi$, the vector $\varphi$ is constant on each orbit of the stabiliser of $v$ in $\operatorname{Aut}(\Fp(p))$.
	
	Since the vertices of $Q_0$ are true twins (having identical closed neighbourhoods), every permutation of $Q_0$ or of $Q_1$ is an automorphism of $\Fp(p)$, as are the translations from $H$. Consequently, the stabiliser of a clique vertex $q\in Q_0$ partitions $V$ into five orbits, and the stabiliser of an outer vertex $x=x_{0,c,t}\in X_0$ partitions $V$ into eight orbits.
	
	\begin{lemma}\label{lem:Fentries}
		Let $p\ge4$ be even, $h=p/2$ and $D=8p^2(p-2)$. For a clique vertex $q\in Q_0$ and an outer vertex $x=x_{0,c,t}\in X_0\cap Y_c$, the entries of $\Lp e_q$ and $\Lp e_x$ on the orbits of the stabilisers of $q$ and of $x$, and the resulting effective resistances, are as follows (E: edge, N: non-edge):
		\begin{center}
			\begin{adjustbox}{max width=\textwidth}
				\renewcommand{\arraystretch}{1.35}
				\begin{tabular}{@{}lccc@{}}
					\toprule
					orbit of $u$ (size) & & $D\,\Lp_{uq}$ & $\rres{u}{q}$\\
					\midrule
					$\{q\}$ ($1$) & & $4p^2-3p-2$ & $0$\\
					$Q_0\setminus\{q\}$ ($p-1$) & E & $5p-2$ & $\frac{1}{p}$\\
					$X_0$ ($p$) & E & $3p+2$ & $\frac{2p+1}{2p^2}$\\
					$X_1$ ($p$) & N & $-(5p-2)$ & $\frac{2p^2+p-2}{2p^2(p-2)}$\\
					$Q_1$ ($p$) & N & $-(7p-6)$ & $\frac{(p-1)(p+2)}{p^2(p-2)}$\\
					\bottomrule
				\end{tabular}\qquad
				\begin{tabular}{@{}lccc@{}}
					\toprule
					orbit of $u$ (size) & & $D\,\Lp_{ux}$ & $\rres{u}{x}$\\
					\midrule
					$\{x\}$ ($1$) & & $4p^2-3p-2$ & $0$\\
					$Q_0$ ($p$) & E & $3p+2$ & $\frac{2p+1}{2p^2}$\\
					$(X_0\cap Y_c)\setminus\{x\}$ ($h-1$) & E & $\frac{3p^2-3p+2}{p-1}$ & $\frac{2p-1}{2p(p-1)}$\\
					$(X_1\cap Y_c)\setminus\{x_{1,c,t}\}$ ($h-1$) & E & $-\frac{p^2+3p-2}{p-1}$ & $\frac{2p^2-3p+2}{2p(p-1)(p-2)}$\\
					$\{x_{1,c,t}\}$ ($1$) & N & $-(3p+2)$ & $\frac{1}{p-2}$\\
					$X_0\cap Y_{1-c}$ ($h$) & N & $-\frac{p^2-9p+6}{p-1}$ & $\frac{2p^2+p-2}{2p^2(p-1)}$\\
					$X_1\cap Y_{1-c}$ ($h$) & N & $-\frac{5p^2-9p+6}{p-1}$ & $\frac{2p^3-p^2-4p+4}{2p^2(p-1)(p-2)}$\\
					$Q_1$ ($p$) & N & $-(5p-2)$ & $\frac{2p^2+p-2}{2p^2(p-2)}$\\
					\bottomrule
				\end{tabular}
			\end{adjustbox}
		\end{center}
		In particular, $\Lp_{qq}=\Lp_{xx}=(4p^2-3p-2)/D$ for every vertex.
	\end{lemma}
	
	\begin{proof}
		For $v=q\in Q_0$, denote the values of $\varphi=\Lp e_q$ on the five orbits $\{q\}$, $Q_0\setminus\{q\}$, $X_0$, $X_1$ and $Q_1$ by $a,\kappa,\alpha,\varepsilon$ and $\zeta$. A vertex in $Q_0\setminus\{q\}$ is adjacent to $q$, to $p-2$ other vertices in $Q_0$ and to all $p$ vertices of $X_0$; a vertex in $X_0$ is adjacent to $q$, to the $p-1$ vertices of $Q_0\setminus\{q\}$, to $h-1$ vertices in $X_0$ and to $h-1$ vertices in $X_1$; and similarly on block $1$. Writing out $(L\varphi)_w=\delta_{qw}-1/(4p)$ on each orbit together with $\one^{\top}\varphi=0$ gives
		\begin{gather*}
			(2p-1)a-(p-1)\kappa-p\alpha=1-\tfrac1{4p},\qquad
			(2p-1)\kappa-a-(p-2)\kappa-p\alpha=-\tfrac1{4p},\\
			(2p-2)\alpha-a-(p-1)\kappa-(h-1)\alpha-(h-1)\varepsilon=-\tfrac1{4p},\\
			(2p-2)\varepsilon-p\zeta-(h-1)\varepsilon-(h-1)\alpha=-\tfrac1{4p},\\
			(2p-1)\zeta-(p-1)\zeta-p\varepsilon=-\tfrac1{4p},\qquad
			a+(p-1)\kappa+p\alpha+p\varepsilon+p\zeta=0.
		\end{gather*}
		Subtracting the second equation from the first immediately gives $a-\kappa=1/(2p)$, and solving the remaining linear equations yields the five values in the left-hand table of Lemma~\ref{lem:Fentries}.
		
		For $v=x=x_{0,c,t}\in X_0\cap Y_c$, the eight stabiliser orbits are listed in the right-hand table of Lemma~\ref{lem:Fentries}. Counting neighbours between orbits directly from Definition~\ref{def:family} (for instance, a vertex $x_{0,c,t'}\in(X_0\cap Y_c)\setminus\{x\}$ is adjacent to all $p$ vertices of $Q_0$, to $x$, to the remaining $h-2$ vertices of $X_0\cap Y_c$, to the $h-2$ vertices $x_{1,c,t''}$ with $t''\notin\{t,t'\}$, and to the partner $x_{1,c,t}$) produces eight linear equations from $L\varphi=e_x-\one/(4p)$ along with $\one^{\top}\varphi=0$, whose unique solution is given in the right-hand table. Finally, since $\Lp_{uu}=\Lp_{vv}=(4p^2-3p-2)/D$ for all $u,v\in V$, substituting $\Lp_{uv}$ into \eqref{eq:resistance} gives the effective resistances in the last column.
	\end{proof}
	
	\subsection{Proof of Theorem~\ref{thm:existence}}
	
	\begin{proof}[Proof of Theorem~\ref{thm:existence}]
		Let $G=\Fp(p)$ with $p\ge4$ even. As noted after Definition~\ref{def:family}, $G$ has $n=4p$ vertices, $2p$ of degree $2p-1$ (in $Q_0\cup Q_1$) and $2p$ of degree $2p-2$ (in $X_0\cup X_1$). By Lemma~\ref{lem:Fbicayley}, $G$ is a bi-Cayley graph over $H$, and every class of nontrivial characters in the table of that lemma has $\Theta>0$ because $p\ge4$, so $G$ is connected by Lemma~\ref{lem:bicayley}. Multiplying the last column of the table by the number of characters in each class and summing gives $-1/p+1/(2p)+(p-1)/(2p(p-1))=0$, so \eqref{eq:criterion} holds and $G$ is resistance-regular by Lemma~\ref{lem:bicayley}. Substituting $|H|=2p$, $|S_{12}|=p$ and the values of $M_{22}/\Theta$ from the table into \eqref{eq:a1} gives $\Lp_{vv}=1/(8p^2)+\bigl(1+1/(p-2)\bigr)/(2p)=(4p^2-3p-2)/(8p^2(p-2))$. By \eqref{eq:rowsum} and $\tr\Lp=4p\,\Lp_{vv}$, every vertex has resistance degree $R_v=(4p^2-3p-2)/(p(p-2))$, and $\Kf(G)=2(4p^2-3p-2)/(p-2)$.
		
		The Laplacian spectrum consists of the eigenvalues of the blocks $\widehat L_\chi$. The trivial character gives $\widehat L_1=p\bigl(\begin{smallmatrix}1&-1\\-1&1\end{smallmatrix}\bigr)$, with eigenvalues $0$ and $2p$. For $(\alpha,\gamma,k)=(0,1,0)$, the block $\bigl(\begin{smallmatrix}2p-2&-p\\-p&p\end{smallmatrix}\bigr)$ has trace $3p-2$ and determinant $p(p-2)$, so its eigenvalues are $\rho_\pm=\tfrac12\bigl(3p-2\pm s\bigr)$ with $s=\sqrt{5p^2-4p+4}$. The remaining four classes have $\chi(S_{12})=0$, so their blocks are diagonal with diagonal entries $(p,2p)$, $(2p-2,2p)$, $(2p,2p)$ (occurring $p-2$ times) and $(2p-2,2p)$ (occurring $p-2$ times). Hence the Laplacian spectrum of $G$ is
		\begin{equation}\label{eq:Fspectrum}
			0^{(1)},\quad p^{(1)},\quad \rho_-^{(1)},\quad \rho_+^{(1)},\quad (2p-2)^{(p-1)},\quad (2p)^{(3p-3)},
		\end{equation}
		where exponents denote multiplicities. Since $s^2-(2p)^2=(p-2)^2>0$, we have $s>2p$, so $\rho_+>(5p-2)/2>2p$ and $\rho_-<(p-2)/2<p$, while $\rho_->0$ because $\rho_+\rho_-=p(p-2)>0$. Together with $p<2p-2<2p$, this gives $0<\rho_-<p<2p-2<2p<\rho_+$, so $G$ has exactly five distinct nonzero Laplacian eigenvalues.
		
		Finally, the translation by $(0,1,0)\in H$ exchanges the two blocks, so every vertex is equivalent under an automorphism to a vertex of $Q_0$ or of $X_0$. By Lemma~\ref{lem:Fentries}, every effective resistance between distinct vertices of $G$ is therefore one of the rational functions of $p$ in the last columns of the two tables of that lemma. The values $(2p+1)/(2p^2)$ and $(2p^2+p-2)/(2p^2(p-2))$ occur in both tables, so there are four of these functions on edges and five on non-edges, and the nine functions are pairwise distinct for every even $p\ge4$.
	\end{proof}
	
	For $p=4,6,8,10$, Theorem~\ref{thm:existence} gives $\Lp_{vv}=25/128$, $31/288$, $115/1536$ and $23/400$, and graphs~1, 4 and~13 of Table~\ref{tab:examples} are $\Fp(4)$, $\Fp(6)$ and $\Fp(8)$. The tables of Lemma~\ref{lem:Fentries} also give a second, combinatorial proof of resistance-regularity, and they explain how the two degree classes balance.
	
	\begin{corollary}\label{cor:resistances}
		In $\Fp(p)$, $\sum_{u\sim q}\rres{u}{q}=\sum_{u\sim x}\rres{u}{x}=2-2/(4p)$ for $q\in Q_b$ and $x\in X_b$; by Theorem~\ref{thm:zwb}(iii), this gives a second proof that $\Fp(p)$ is resistance-regular.
	\end{corollary}
	
	The two sums in Corollary~\ref{cor:resistances} follow from Lemma~\ref{lem:Fentries} and show how $\Fp(p)$ achieves constant resistance degree despite having two different vertex degrees. Every vertex, whether in $Q_b$ or in $X_b$, is incident to $p$ cross-edges between $Q_b$ and $X_b$, and each such edge has resistance $(2p+1)/(2p^2)$, contributing $(2p+1)/(2p)$ to both sums. The remaining incident edges are the $p-1$ clique edges at $q\in Q_b$ and the $p-2=2(h-1)$ cocktail-party edges at $x\in X_b$. Their contributions match because of the identity
	\begin{equation}\label{eq:balance}
		(p-1)\cdot\frac1p=(h-1)\Bigl(\frac{2p-1}{2p(p-1)}+\frac{2p^2-3p+2}{2p(p-1)(p-2)}\Bigr)=\frac{p-1}{p}.
	\end{equation}
	In view of \eqref{eq:kirchhoff}, equation~\eqref{eq:balance} states that in a uniformly random spanning tree of $\Fp(p)$, the expected number of incident clique edges at a vertex of $Q_b$ equals the expected number of incident $Y_c$-edges at a vertex of $X_b$: the $p-2$ edges of $Y_c$ at $x$ are fewer in number than the $p-1$ clique edges at $q$, but each carries a slightly higher spanning-tree probability because the removed matching in $Y_c$ forces current between the two blocks to pass through them. Note also that the value $\rres{q}{q'}=1/p$ on a clique edge has a transparent electrical reason: any two vertices $q,q'\in Q_b$ are true twins, so when unit current flows from $q$ to $q'$ every other vertex stays at the midpoint potential, and the direct edge $qq'$ acts in parallel with $2p-2$ two-step paths through common neighbours, giving resistance $1/(1+(2p-2)/2)=1/p$.
	
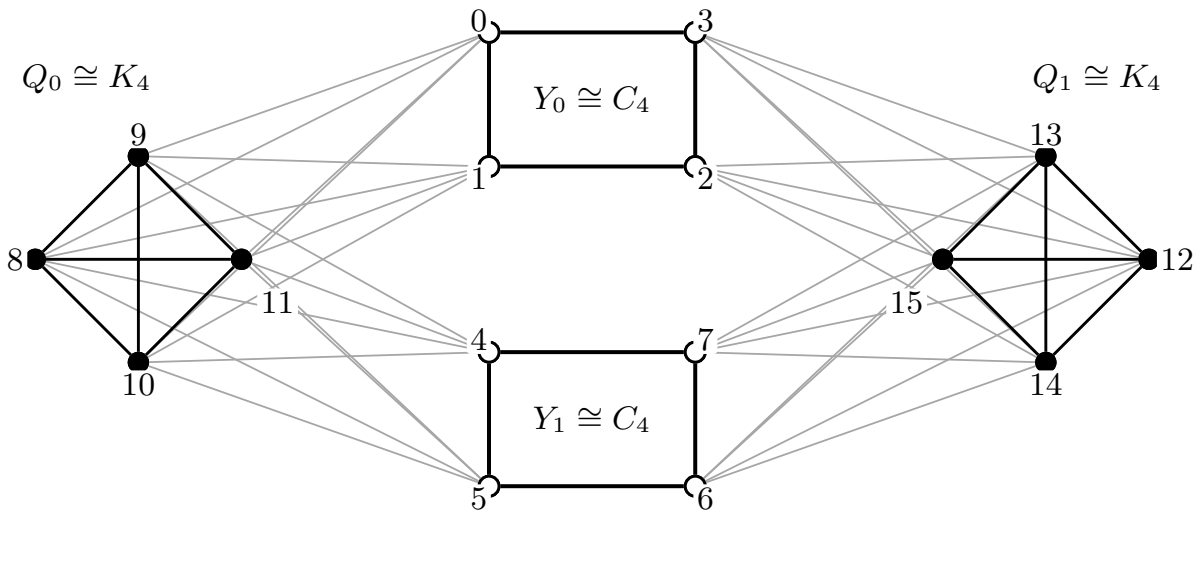
\begin{figure}[!htbp]
		\centering
		\begin{adjustbox}{width=0.5\textwidth}
			\begin{tikzpicture}[scale=1.0,
				qv/.style={circle, draw, fill=black, inner sep=0pt, minimum size=5.5pt},
				xv/.style={circle, draw, fill=white, inner sep=0pt, minimum size=5.5pt, thick},
				jn/.style={gray!70, line width=0.5pt},
				cyc/.style={line width=1.1pt},
				clq/.style={line width=0.8pt},
				lab/.style={font=\footnotesize, fill=white, inner sep=1pt, rounded corners=1pt}]
				\node[qv] (q8)  at (-5.4, 0)   {}; \node[qv] (q9)  at (-4.4, 1)  {};
				\node[qv] (q10) at (-4.4,-1)   {}; \node[qv] (q11) at (-3.4, 0)  {};
				\node[qv] (q12) at ( 5.4, 0)   {}; \node[qv] (q13) at ( 4.4, 1)  {};
				\node[qv] (q14) at ( 4.4,-1)   {}; \node[qv] (q15) at ( 3.4, 0)  {};
				\node[xv] (v0) at (-1.0, 2.2) {}; \node[xv] (v1) at (-1.0, 0.9) {};
				\node[xv] (v2) at ( 1.0, 0.9) {}; \node[xv] (v3) at ( 1.0, 2.2) {};
				\node[xv] (v4) at (-1.0,-0.9) {}; \node[xv] (v5) at (-1.0,-2.2) {};
				\node[xv] (v6) at ( 1.0,-2.2) {}; \node[xv] (v7) at ( 1.0,-0.9) {};
				\foreach \q in {q8,q9,q10,q11} \foreach \v in {v0,v1,v4,v5} \draw[jn] (\q) -- (\v);
				\foreach \q in {q12,q13,q14,q15} \foreach \v in {v3,v2,v7,v6} \draw[jn] (\q) -- (\v);
				\foreach \a/\b in {q8/q9,q8/q10,q8/q11,q9/q10,q9/q11,q10/q11,q12/q13,q12/q14,q12/q15,q13/q14,q13/q15,q14/q15} \draw[clq] (\a) -- (\b);
				\draw[cyc] (v0) -- (v1) -- (v2) -- (v3) -- (v0);
				\draw[cyc] (v4) -- (v5) -- (v6) -- (v7) -- (v4);
				\foreach \q in {q8,q9,q10,q11,q12,q13,q14,q15} \node[qv] at (\q) {};
				\foreach \v in {v0,v1,v2,v3,v4,v5,v6,v7} \node[xv] at (\v) {};
				\node[lab, left=2pt] at (q8) {$8$}; \node[lab, above=2pt] at (q9) {$9$}; \node[lab, below=2pt] at (q10) {$10$}; \node[lab] at (-3.05, -0.42) {$11$};
				\node[lab, right=2pt] at (q12) {$12$}; \node[lab, above=2pt] at (q13) {$13$}; \node[lab, below=2pt] at (q14) {$14$}; \node[lab] at (3.05, -0.42) {$15$};
				\node[lab, above left=-1pt] at (v0) {$0$}; \node[lab, below left=-1pt] at (v1) {$1$}; \node[lab, below right=-1pt] at (v2) {$2$}; \node[lab, above right=-1pt] at (v3) {$3$};
				\node[lab, above left=-1pt] at (v4) {$4$}; \node[lab, below left=-1pt] at (v5) {$5$}; \node[lab, below right=-1pt] at (v6) {$6$}; \node[lab, above right=-1pt] at (v7) {$7$};
				\node[lab, fill=none] at (-4.9, 1.75) {$Q_0\cong K_4$}; \node[lab, fill=none] at (4.9, 1.75) {$Q_1\cong K_4$};
				\node[lab, fill=none] at (0, 1.55) {$Y_0\cong C_4$}; \node[lab, fill=none] at (0,-1.55) {$Y_1\cong C_4$};
			\end{tikzpicture}
		\end{adjustbox}
		\caption{The graph $\Fp(4)$. Filled vertices form the two cliques $Q_0$ and $Q_1$ (degree $7$), hollow vertices the outer sets $X_0=\{0,1,4,5\}$ and $X_1=\{2,3,6,7\}$ (degree $6$). Thin grey lines are the complete joins $Q_0$--$X_0$ and $Q_1$--$X_1$; the heavy squares are the $4$-cycles induced by $Y_0$ and $Y_1$.}
		\label{fig:f4}
	\end{figure}

	\begin{example}[The sixteen-vertex graph $\Fp(4)$]\label{ex:16}
		The smallest member $\Fp(4)$ has $16$ vertices and $52$ edges, eight vertices of degree $6$ and eight of degree $7$. It is drawn in Figure~\ref{fig:f4} with $Q_0=\{8,9,10,11\}$, $Q_1=\{12,13,14,15\}$, $X_0=\{0,1,4,5\}$ and $X_1=\{2,3,6,7\}$, where $Y_0=\{0,1,2,3\}$ and $Y_1=\{4,5,6,7\}$ induce $4$-cycles. Here $\Lp_{vv}=25/128$, $R_v=25/4$ and $\Kf(\Fp(4))=50$; the Laplacian eigenvalues are $0$, $4$, $6^{(3)}$, $8^{(9)}$ and $5\pm\sqrt{17}$, and $\Fp(4)$ has $\tau=2^{31}\cdot3^{3}$ spanning trees, in accordance with Proposition~\ref{prop:necessary}. By Lemma~\ref{lem:Fentries}, the resistances on edges are $1/4$, $9/32$, $7/24$ and $11/24$, and those on non-edges are $17/48$, $1/2$, $25/48$, $17/32$ and $9/16$. Thus $\rres{u}{v}$ determines whether $uv$ is an edge, but no threshold on $\rres{u}{v}$ separates edges from non-edges. The automorphism group of $\Fp(4)$ has order $9216$, and its orbits are the two degree classes.
	\end{example}
	
	\paragraph{A four-parameter family.} The graph $\Fp(p)$ is a special case of a larger construction. For integers $s,a,k\ge1$ and $h\ge2$, take $s$ disjoint cliques $Q_b\cong K_a$ ($b<s$), join each $Q_b$ completely to a set $X_b$ of $kh$ new vertices $x_{b,c,t}$ ($c<k$, $t<h$), and join $x_{b,c,t}$ and $x_{b',c',t'}$ whenever $c=c'$ and $t\ne t'$. Writing $\mathcal G(s,a,k,h)$ for this graph, we have $\Fp(p)=\mathcal G(2,p,2,p/2)$. In Appendix~\ref{app:fourparam} we compute the Laplacian spectrum of $\mathcal G(s,a,k,h)$ and show that it is resistance-regular if and only if a single polynomial $N(s,a,k,h)$ vanishes. The integer solutions include a two-parameter family in which the two degrees differ by $2s-3$; since the constructions of Section~\ref{sec:far} give larger differences, we do not pursue this family in the main text.
	
	\subsection{How sharp Theorem~\ref{thm:forced} is}\label{sec:sharp}
	
	The graphs $\Fp(p)$ have $m=5$ distinct nonzero Laplacian eigenvalues. Graphs~14 and~15 of Table~\ref{tab:examples} have one fewer.
	
	\begin{proposition}\label{prop:sharp}
		There exist non-regular resistance-regular graphs whose Laplacian has $m=4$ distinct nonzero eigenvalues.
	\end{proposition}
	
	\begin{proof}
		Let $G$ be graph~14 (or graph~15) of Table~\ref{tab:examples}, defined as a bi-Cayley graph over $\Z_4\times\Z_4$ by the connection sets in Appendix~\ref{app:examples}. It has $32$ vertices of degrees $14$ and $15$, and exact rational inversion of $L+J/32$ confirms that $\Lp_{vv}=37/512$ for every vertex $v$. Its Laplacian characteristic polynomial is
		\[
		x\,(x-12)^{6}\,(x-16)^{23}\,(x^2-24x+64),
		\]
		so its nonzero Laplacian eigenvalues are $12$, $16$ and $12\pm4\sqrt5$, which are four distinct values.
	\end{proof}
	
	Hence, for $m=4$ and $m=5$, the conditions in Theorem~\ref{thm:forced}(i) cannot be reduced to the first one, even if only regularity, rather than walk-regularity, is required. For $m=3$, Theorem~\ref{thm:forced}(ii) leaves room for two distinct degrees, and we do not know whether a non-regular resistance-regular graph with $m=3$ exists. Whether the first two conditions in Theorem~\ref{thm:forced}(i) always force regularity is discussed in Section~\ref{sec:moment}.

	\section{Modules and blow-ups}\label{sec:modules}
	
	The proofs in Section~\ref{sec:construction} used one structural feature of $\Fp(p)$ repeatedly: the vertices of each clique $Q_b$ are true twins, so every permutation of $Q_b$ is an automorphism and the potentials inside $Q_b$ are easy to control. In this section we turn this observation into a general method, which we use in Section~\ref{sec:far} to prove Theorem~\ref{thm:far}. We show that the diagonal entry $\Lp_{vv}$ at a vertex $v$ of a module splits into a part that depends only on the module itself and a part that depends only on how the module is attached to the rest of the graph (Lemma~\ref{lem:module}). For graphs built entirely from cliques and independent sets, resistance-regularity then becomes a system of rational equations in the sizes of these sets.
	
	\begin{definition}[module, blow-up]\label{def:module}
		A set $M\subseteq V$ is a \emph{module} of $G$ if every vertex outside $M$ is adjacent either to all vertices of $M$ or to none of them. Given a connected graph $H$ on the vertex set $\{1,\dots,r\}$ and graphs $W_1,\dots,W_r$, the \emph{blow-up} $H[W_1,\dots,W_r]$ is obtained from the disjoint union of $W_1,\dots,W_r$ by joining every vertex of $W_i$ to every vertex of $W_j$ whenever $ij$ is an edge of $H$.
	\end{definition}
	
	The vertex sets $M_i=V(W_i)$ of a blow-up are modules, which we call the modules of the blow-up. We write $K_m$ for the complete graph and $\overline{K_m}$ for the edgeless graph on $m$ vertices; the vertices of a module inducing $K_m$ are true twins (they have the same closed neighbourhood), and those of a module inducing $\overline{K_m}$ are false twins (they have the same open neighbourhood).
	
	\begin{lemma}\label{lem:module}
		Let $G$ be a connected graph on $n$ vertices, let $M\ne V$ be a module of $G$ with $|M|=m$ and $W=G[M]$, and let $c\ge1$ be the number of vertices outside $M$ that are adjacent to $M$. Write $L_W$ for the Laplacian of $W$ and $\one_M$ for the indicator vector of $M$. Then for every $v\in M$,
		\begin{align}
			\Lp_{vv}&=\bigl[(L_W+cI)^{-1}\bigr]_{vv}-\frac{1}{cm}+\frac{\one_M^{\top}\Lp\one_M}{m^2},\label{eq:module1}\\
			\bigl((\Lp)^2\bigr)_{vv}&=\Bigl\|(L_W+cI)^{-1}\Bigl(e_v-\frac{\one_M}{m}\Bigr)\Bigr\|^2+\frac{\|\Lp\one_M\|^2}{m^2}.\label{eq:module2}
		\end{align}
		If $W$ is walk-regular, then $\Lp_{vv}$ and $\bigl((\Lp)^2\bigr)_{vv}$ take the same value at all vertices of $M$, and
		\[
		\bigl[(L_W+cI)^{-1}\bigr]_{vv}-\frac1{cm}=\frac1m\sum_{\nu}\frac{1}{\nu+c},
		\]
		where $\nu$ runs over all Laplacian eigenvalues of $W$, with multiplicity, except for one eigenvalue $0$.
	\end{lemma}
	
	\begin{proof}
		Let $U$ be the space of vectors supported on $M$ whose entries sum to zero. For $z\in U$ and a vertex $w\notin M$, $(Lz)_w=-\sum_{u\sim w}z_u$ is either $0$ or $-\sum_{u\in M}z_u=0$, because $w$ is adjacent to all of $M$ or to none of it; and for $w\in M$, $(Lz)_w=\bigl((L_W+cI)z\bigr)_w$, because $w$ has exactly $c$ neighbours outside $M$, at which $z$ vanishes. Hence $U$ is invariant under $L$, and $L$ acts on $U$ as $L_W+cI$, which is invertible because $c\ge1$. Since $L$ is symmetric and $U\subseteq\one^\perp$, the space $U^\perp$ is also invariant, so $\Lp$ maps $U$ to itself, acting there as $(L_W+cI)^{-1}$, and maps $U^\perp$ to itself. Write $e_v=z+\one_M/m$ with $z=e_v-\one_M/m\in U$, and note that $\one_M\in U^\perp$. Then $\Lp_{vv}=z^\top\Lp z+\one_M^\top\Lp\one_M/m^2$ and $\bigl((\Lp)^2\bigr)_{vv}=\|\Lp e_v\|^2=\|\Lp z\|^2+\|\Lp\one_M\|^2/m^2$, the cross terms vanishing because $\Lp z\in U$ and $\Lp\one_M\in U^\perp$. Finally, $L_W\one_M=0$ gives $(L_W+cI)^{-1}\one_M=\one_M/c$, so that $z^\top(L_W+cI)^{-1}z=\bigl[(L_W+cI)^{-1}\bigr]_{vv}-2/(cm)+1/(cm)$, which proves \eqref{eq:module1}; \eqref{eq:module2} is $\|\Lp z\|^2=\|(L_W+cI)^{-1}z\|^2$.
		
		If $W$ is walk-regular of degree $d$, then $L_W=dI-A_W$, so $(L_W+cI)^{-1}$ and $(L_W+cI)^{-2}$ are polynomials in $A_W$ and have constant diagonals; the remaining terms in \eqref{eq:module1} and \eqref{eq:module2} do not depend on $v\in M$. The constant diagonal of $(L_W+cI)^{-1}$ equals $\frac1m\tr(L_W+cI)^{-1}=\frac1m\sum_{\nu'}\frac1{\nu'+c}$, the sum running over all Laplacian eigenvalues $\nu'$ of $W$; removing one eigenvalue $0$ cancels the term $1/(cm)$.
	\end{proof}
	
	For a clique module $K_m$ the sum in Lemma~\ref{lem:module} consists of $m-1$ terms $1/(m+c)$, and for an independent module $\overline{K_m}$ of $m-1$ terms $1/c$. We now combine the lemma over all modules of a blow-up. Let $H$ be a connected graph on $\{1,\dots,r\}$ with $r\ge2$, let $m_1,\dots,m_r\ge1$, let each $W_i$ be $K_{m_i}$ or $\overline{K_{m_i}}$, and let $G=H[W_1,\dots,W_r]$, so that $n=m_1+\dots+m_r$. For each $i$ put
	\[
	c_i=\sum_{j\sim i}m_j,\qquad \lambda_i=\begin{cases}c_i+m_i&\text{if }W_i=K_{m_i},\\ c_i&\text{if }W_i=\overline{K_{m_i}},\end{cases}
	\]
	so that a vertex of $M_i$ has degree $\lambda_i-1$ if $M_i$ is a clique and $\lambda_i$ if it is independent.
	
	\begin{proposition}\label{prop:blowup}
		For each $i$, the linear system
		\begin{equation}\label{eq:quotient}
			\sum_{k\sim j}m_k\,(y_j-y_k)=\delta_{ij}-\frac{m_i}{n}\quad(1\le j\le r),\qquad \sum_{j=1}^r m_jy_j=0,
		\end{equation}
		has a unique solution $y^{(i)}\in\R^r$, and for every vertex $v\in M_i$,
		\[
		\Lp_{vv}=\frac{m_i-1}{m_i\lambda_i}+\frac{y^{(i)}_i}{m_i},\qquad
		\bigl((\Lp)^2\bigr)_{vv}=\frac{m_i-1}{m_i\lambda_i^2}+\frac{1}{m_i^2}\sum_{j=1}^r m_j\bigl(y^{(i)}_j\bigr)^2.
		\]
		In particular, $G$ is resistance-regular if and only if $(m_i-1)/(m_i\lambda_i)+y^{(i)}_i/m_i$ takes the same value for all $i$.
	\end{proposition}
	
	\begin{proof}
		A vector that is constant on every module, with value $y_j$ on $M_j$, is mapped by $L$ to the vector whose value on $M_j$ is $\sum_{k\sim j}m_k(y_j-y_k)$, because the edges inside a module do not contribute. Hence the space of such vectors is invariant under $L$, and so is its intersection with $\one^\perp$, on which $L$ is invertible. The vector $\one_{M_i}-(m_i/n)\one$ lies in this intersection, so $\Lp\one_{M_i}$ is constant on every module; writing $y_j$ for its value on $M_j$, the equations $L\Lp\one_{M_i}=\one_{M_i}-(m_i/n)\one$ and $\one^\top\Lp\one_{M_i}=0$ are exactly \eqref{eq:quotient}. The solution is unique, because the difference of two solutions corresponds to a vector in the kernel of $L$, which is constant and orthogonal to $\one$, hence zero. Now $\one_{M_i}^\top\Lp\one_{M_i}=m_iy^{(i)}_i$ and $\|\Lp\one_{M_i}\|^2=\sum_jm_j(y^{(i)}_j)^2$, and the first terms of \eqref{eq:module1} and \eqref{eq:module2} equal $(m_i-1)/(m_i\lambda_i)$ and $(m_i-1)/(m_i\lambda_i^2)$, because $(L_{W_i}+c_iI)^{-1}$ acts on the vectors in $\R^{M_i}$ with zero sum as multiplication by $1/\lambda_i$.
	\end{proof}
	
	Effective resistances in several classes of graph blow-ups have been studied previously. Azimi, Bapat and Farrokhi \cite{AzimiBapatFarrokhi2021} derived resistance formulas for blow-ups of trees in which every vertex is replaced by a copy of a fixed clique or a fixed independent set. More recently, Liu, Li, You, Hua and Chen \cite{LiuLiYouHuaChen2025} treated generalised blow-ups of an arbitrary base graph $H$ whose modules are disjoint unions of isolated vertices and copies of a single clique $K_t$ (with $t$ common to all modules), reducing their effective resistances to those of a vertex-weighted quotient on $V(H)$. Because $\Lp_{vv}=R_v/n-\Kf(G)/n^2$ by \eqref{eq:rowsum}, pairwise resistance formulas determine the diagonal of $\Lp$ in principle for those families. Proposition~\ref{prop:blowup} accommodates cliques and independent sets of arbitrary, module-dependent sizes $m_i$ and isolates the diagonals of both $\Lp$ and $(\Lp)^2$ directly in the form required to test resistance-regularity.
	
	When $H$ is a tree, the system \eqref{eq:quotient} can be solved by inspection: the current between two adjacent modules is determined by the total demand on one side, as in a series circuit. This gives a criterion that can be checked by hand.
	
	\begin{corollary}\label{cor:tree}
		Let $H$ be a tree on $\{1,\dots,r\}$, let $m_1,\dots,m_r\ge1$, let $W_1,\dots,W_r$ be walk-regular graphs with $|V(W_i)|=m_i$, and let $G=H[W_1,\dots,W_r]$. For each $i$ put $T_i=\frac1{m_i}\sum_\nu\frac{1}{\nu+c_i}$ as in Lemma~\ref{lem:module}, with $c_i=\sum_{j\sim i}m_j$, and for an edge $jk$ of $H$ and a vertex $i$ of $H$ let $A^{(i)}_{jk}$ be the total size of the modules in the component of $H-jk$ that does not contain $i$. Then for every vertex $v\in M_i$,
		\[
		\Lp_{vv}=T_i+\frac{1}{n^2}\sum_{jk\in E(H)}\frac{\bigl(A^{(i)}_{jk}\bigr)^2}{m_jm_k}.
		\]
		Consequently $G$ is resistance-regular if and only if, for every edge $ij$ of $H$,
		\begin{equation}\label{eq:tree}
			T_j-T_i=\frac{n-2A_{ij}}{n\,m_im_j},
		\end{equation}
		where $A_{ij}$ is the total size of the modules on the side of $i$ in $H-ij$.
	\end{corollary}
	
	\begin{proof}
		By Lemma~\ref{lem:module}, $\Lp_{vv}=T_i+\one_{M_i}^\top\Lp\one_{M_i}/m_i^2$. Set $\varphi=\Lp\one_{M_i}$. By the invariance argument in the proof of Proposition~\ref{prop:blowup} (which does not depend on the internal structure of the modules), $\varphi$ is constant on every module $M_j$, and hence only edges between distinct modules contribute to the quadratic form $\one_{M_i}^\top\Lp\one_{M_i}=\varphi^\top L\varphi=\sum_{ab\in E}(\varphi_a-\varphi_b)^2$. Let $jk$ be an edge of $H$ and $S$ the set of vertices of $G$ lying in modules of the component of $H-jk$ not containing $i$. Summing the equations $(L\varphi)_a=-m_i/n$ over $a\in S$, the edges inside $S$ cancel and only the $m_jm_k$ edges between $M_j$ and $M_k$ remain, each carrying the same potential difference; hence $m_jm_k|\varphi_j-\varphi_k|=m_iA^{(i)}_{jk}/n$, where $\varphi_j$ denotes the value of $\varphi$ on $M_j$. The energy of $\varphi$ is therefore $\sum_{jk}m_jm_k(\varphi_j-\varphi_k)^2=(m_i^2/n^2)\sum_{jk}(A^{(i)}_{jk})^2/(m_jm_k)$, which gives the formula for $\Lp_{vv}$. If $i$ and $j$ are adjacent in $H$, then every edge of $H$ other than $ij$ has the same far side as seen from $i$ and from $j$, while for the edge $ij$ itself the far side has size $n-A_{ij}$ as seen from $i$ and $A_{ij}$ as seen from $j$. Hence $\Lp_{vv}-\Lp_{ww}=T_i-T_j+\bigl((n-A_{ij})^2-A_{ij}^2\bigr)/(n^2m_im_j)$ for $v\in M_i$ and $w\in M_j$, which vanishes exactly when \eqref{eq:tree} holds. Since $H$ is connected, $\Lp$ has constant diagonal if and only if this holds for every edge.
	\end{proof}
	
	\paragraph{Graphs with two consecutive degrees.} A different substitution turns any example with two consecutive degrees into an infinite sequence of examples.
	
	\begin{proposition}\label{prop:gap1}
		Let $G$ be a connected resistance-regular graph on $n$ vertices whose degrees are $d$ and $d+1$, and let $m\ge1$. Let $G^{(m)}$ be the blow-up of $G$ in which every vertex of degree $d$ is replaced by a clique $K_m$ and every vertex of degree $d+1$ by an independent set $\overline{K_m}$. Then $G^{(m)}$ is a connected resistance-regular graph on $mn$ vertices with degrees $m(d+1)-1$ and $m(d+1)$, and for a vertex $w$ in the module of a vertex $v$ of $G$,
		\[
		\Lp_{ww}(G^{(m)})=\frac{m-1}{m^2(d+1)}+\frac{\Lp_{vv}(G)}{m^2},\qquad
		\bigl((\Lp)^2\bigr)_{ww}(G^{(m)})=\frac{m-1}{m^3(d+1)^2}+\frac{\bigl((\Lp)^2\bigr)_{vv}(G)}{m^3}.
		\]
	\end{proposition}
	
	\begin{proof}
		All modules have size $m$. A module replacing a vertex $v$ has $c_v=md_v$ outside neighbours, so $\lambda_v=md+m=m(d+1)$ if $d_v=d$ (clique) and $\lambda_v=m(d+1)$ if $d_v=d+1$ (independent set); thus $\lambda_v=m(d+1)$ for every $v$. With all $m_j=m$, the system \eqref{eq:quotient} for the module of $v$ reads $m\,(L_Gy)_u=\delta_{uv}-1/n$ for $u\in V(G)$, together with $\sum_uy_u=0$, so its solution is $y=\frac1m\Lp(G)e_v$. Hence $y_v=\Lp_{vv}(G)/m$ and $\sum_um\,y_u^2=\bigl((\Lp)^2\bigr)_{vv}(G)/m$, and Proposition~\ref{prop:blowup} gives the two formulas; the first is constant because $G$ is resistance-regular. The degrees are $(m-1)+md$ and $m(d+1)$.
	\end{proof}
	
	Applied to the graphs of Table~\ref{tab:examples}, Proposition~\ref{prop:gap1} produces from each of them an infinite sequence of non-regular resistance-regular graphs with two consecutive degrees, and graph~12 of that table is in fact $\Fp(4)^{(2)}$. Several of the other graphs in Table~\ref{tab:examples} are blow-ups as well; for instance, graph~2 is the blow-up $C_9[K_4,\overline{K_2},\overline{K_2},K_4,\overline{K_2},\overline{K_2},K_4,\overline{K_2},\overline{K_2}]$ of the $9$-cycle. By contrast, replacing every vertex of $\Fp(4)$ by modules of a single type, as in the lexicographic products $\Fp(4)[K_2]$ and $\Fp(4)[\overline{K_2}]$, destroys resistance-regularity, and other standard graph products do not produce new examples either.
	
	\begin{remark}[graph products]\label{rem:shape}
		For $G=\Fp(4)$, none of the Cartesian products $G\,\square\,K_2$, $G\,\square\,C_4$, $G\,\square\,K_4$, $G\,\square\,G$ and $G\,\square\,\Fp(6)$, the tensor product $G\times K_2$, the lexicographic products $G[K_2]=G\boxtimes K_2$ and $G[\overline{K_2}]$, the complement and the line graph of $G$ is resistance-regular, and joins are excluded in general by Theorem~\ref{thm:forced}(iii). For a Cartesian product $G\,\square\,H$, the diagonal entry of $\Lp$ at $(v,w)$ is $\sum_{(i,j)\ne(0,0)}e_i(v)f_j(w)/(\mu_i+\nu_j)$, where $(e_i(v))$ and $(f_j(w))$ are the local spectral measures of $G$ and $H$ including the eigenvalue $\mu_0=\nu_0=0$; resistance-regularity of $G$ controls the resolvent $\nu\mapsto\sum_{i\ge1}e_i(v)/(\mu_i+\nu)$ only at $\nu=0$.
	\end{remark}

	\section{How far from regular}\label{sec:far}
	
	Every non-regular resistance-regular graph constructed in Section~\ref{sec:construction} had two consecutive vertex degrees $d$ and $d+1$. Using the module machinery of Section~\ref{sec:modules}, we now show that neither the degree difference $\Delta-\delta$ nor the number of distinct degrees is limited to the patterns seen in Table~\ref{tab:examples}. We first construct a family of blown-up hexagons $C_6(t,t^2)$ whose two degrees differ by $(\sqrt{4n+1}-5)/4$ and path blow-ups with three and four distinct degrees, and then establish the general upper bounds on $\Delta-\delta$ stated in Theorem~\ref{thm:far}(iii). Throughout this section, $\Delta$ and $\delta$ denote the maximum and minimum degrees of $G$, and $R$ also denotes the common resistance degree $R_v$ when $G$ is resistance-regular.
	
	\subsection{The blown-up hexagon}
	
\begin{figure}[!htbp]
		\centering
		\begin{adjustbox}{width=0.45\textwidth}
			\begin{tikzpicture}[scale=1, every node/.style={font=\small}]
				\tikzset{blob/.style={draw, rounded corners=6pt, minimum width=1.6cm, minimum height=0.9cm, align=center, fill=gray!8}}
				\node[blob, fill=gray!25] (M0) at (0,2.2) {$M_0\cong\overline{K_x}$};
				\node[blob] (M1) at (2.6,1.1) {$M_1\cong K_y$};
				\node[blob] (M2) at (2.6,-1.1) {$M_2\cong K_y$};
				\node[blob, fill=gray!25] (M3) at (0,-2.2) {$M_3\cong\overline{K_x}$};
				\node[blob] (M4) at (-2.6,-1.1) {$M_4\cong K_y$};
				\node[blob] (M5) at (-2.6,1.1) {$M_5\cong K_y$};
				\draw[double, thick] (M0) -- (M1);
				\draw[double, thick] (M1) -- (M2);
				\draw[double, thick] (M2) -- (M3);
				\draw[double, thick] (M3) -- (M4);
				\draw[double, thick] (M4) -- (M5);
				\draw[double, thick] (M5) -- (M0);
			\end{tikzpicture}
		\end{adjustbox}
		\caption{The graph $C_6(x,y)$. Shaded modules are independent sets of size $x$, white modules are cliques of size $y$, and double lines indicate complete joins. By Lemma~\ref{lem:c6}, the graph is resistance-regular exactly when $y=x^2$.}
		\label{fig:c6}
	\end{figure}
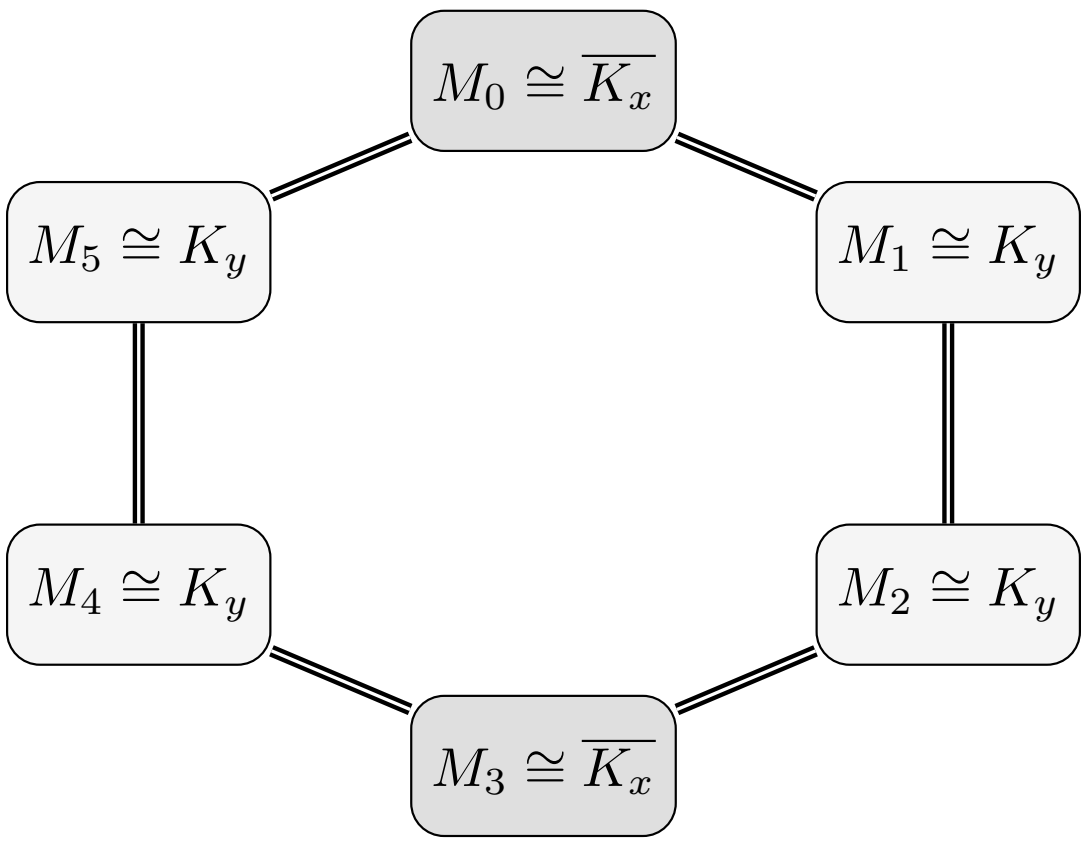

	For integers $x,y\ge1$, let $C_6(x,y)$ be the blow-up of the cycle $C_6$ on $\{0,1,\dots,5\}$ in which $M_0$ and $M_3$ are independent sets of size $x$ and $M_1,M_2,M_4,M_5$ are cliques of size $y$; consecutive modules $M_i$ and $M_{i+1}$ (indices modulo $6$) are completely joined, as in Figure~\ref{fig:c6}. The graph has $n=2x+4y$ vertices; a vertex of $M_0$ or $M_3$ has degree $2y$, and a vertex of a clique module has degree $(y-1)+x+y=x+2y-1$.
	
	\begin{lemma}\label{lem:c6}
		For all integers $x,y\ge1$, a vertex $v$ of an independent module and a vertex $w$ of a clique module of $C_6(x,y)$ satisfy
		\begin{align*}
			\Lp_{vv}&=\frac{4x^2y+x^2+16xy^2+2xy+16y^3-4y^2}{8y^2(x+2y)^2},\\
			\Lp_{ww}&=\frac{x^3+8x^2y^2+2x^2y+16xy^3+8y^3}{8xy^2(x+2y)^2},
		\end{align*}
		and therefore
		\[
		\Lp_{vv}-\Lp_{ww}=\frac{x^2-y}{2xy(x+2y)}.
		\]
		In particular, $C_6(x,y)$ is resistance-regular if and only if $y=x^2$.
	\end{lemma}
	
	\begin{proof}
		The reflection of $C_6$ fixing $0$ and $3$ and the reflection exchanging $0$ with $3$ and $1$ with $2$ preserve the types of the modules, and together with the permutations inside the modules they generate a group of automorphisms that is transitive on $M_0\cup M_3$ and on the union of the four clique modules. It therefore suffices to take $v\in M_0$ and $w\in M_1$. Here $\lambda_0=2y$ and $\lambda_1=x+2y$. Solving \eqref{eq:quotient} for $i=0$ and $i=1$ gives, with $D_0=8y^2(x+2y)^2$ and $D_1=8xy(x+2y)^2$,
		\begin{center}
			\begin{adjustbox}{max width=\textwidth}
				\renewcommand{\arraystretch}{1.2}
				\begin{tabular}{@{}lll@{\hspace{2em}}lll@{}}
					\toprule
					$j$ & $D_0\,y^{(0)}_j$ & $D_1\,y^{(1)}_j$ & $j$ & $D_0\,y^{(0)}_j$ & $D_1\,y^{(1)}_j$\\
					\midrule
					$0$ & $(x+4y)(x^2+2xy+4y^2)$ & $x^2(x+4y)$ & $3$ & $-x(x^2+6xy+12y^2)$ & $-x(x^2+4xy+8y^2)$\\
					$1$ & $x^2(x+4y)$ & $x^3+10x^2y+16xy^2+8y^3$ & $4$ & $-x(x^2+4xy+8y^2)$ & $-(x^3+2x^2y+8xy^2+8y^3)$\\
					$2$ & $-x(x^2+4xy+8y^2)$ & $-(x^3+2x^2y-8y^3)$ & $5$ & $x^2(x+4y)$ & $x^3+2x^2y-8xy^2-8y^3$\\
					\bottomrule
				\end{tabular}
			\end{adjustbox}
		\end{center}
		as one checks by substituting into the six equations of \eqref{eq:quotient}. By Proposition~\ref{prop:blowup},
		\begin{align*}
			\Lp_{vv}&=\frac{x-1}{2xy}+\frac{y^{(0)}_0}{x}=\frac{4x^2y+x^2+16xy^2+2xy+16y^3-4y^2}{8y^2(x+2y)^2},\\
			\Lp_{ww}&=\frac{y-1}{y(x+2y)}+\frac{y^{(1)}_1}{y}=\frac{x^3+8x^2y^2+2x^2y+16xy^3+8y^3}{8xy^2(x+2y)^2},
		\end{align*}
		and subtracting gives the stated difference. All vertices of the independent modules have the same value of $\Lp_{vv}$, and so do all vertices of the clique modules, which proves the last statement.
	\end{proof}
	
	\subsection{Blow-ups of paths}
	
	All examples so far have exactly two distinct degrees. Corollary~\ref{cor:tree} makes it easy to look for blow-ups of paths with more degrees, and the smallest examples we found are the following. For walk-regular graphs $W_1,\dots,W_r$ write $P_r[W_1,\dots,W_r]$ for the blow-up of the path $1-2-\dots-r$.
	
	\begin{lemma}\label{lem:paths}
		The graph
		\[
		\Gamma_{40}=P_8\bigl[K_9,\overline{K_6},\overline{K_3},\overline{K_2},\overline{K_2},\overline{K_3},\overline{K_6},K_9\bigr]
		\]
		is a connected resistance-regular graph on $40$ vertices with exactly four distinct degrees, $5,8,12,14$, carried by $4,6,12,18$ vertices respectively, and $\Lp_{vv}=99/400$. The graph
		\[
		\Gamma_{126}=P_{10}\bigl[K_{33},\overline{K_9},\overline{K_3},K_{3,3,3},\overline{K_9},\overline{K_9},K_{3,3,3},\overline{K_3},\overline{K_9},K_{33}\bigr]
		\]
		is a connected resistance-regular graph on $126$ vertices with exactly three distinct degrees, $18,36,41$, carried by $42,18,66$ vertices respectively, and $\Lp_{vv}=1193/15876$.
	\end{lemma}
	
	\begin{proof}
		Both blow-ups are symmetric under the reversal of the path, so it suffices to check \eqref{eq:tree} for the edges in the first half of the path, the edge in the middle being balanced by symmetry. For $\Gamma_{40}$ the modules $M_1,\dots,M_4$ have $c_i=6,12,8,5$ and $T_i=\frac{8}{135},\frac{5}{72},\frac1{12},\frac1{10}$, and for the edges $\{1,2\},\{2,3\},\{3,4\}$ we have $A_{ij}=9,15,18$, so that
		\[
		\frac{40-18}{40\cdot9\cdot6}=\frac{11}{1080}=T_2-T_1,\quad \frac{40-30}{40\cdot6\cdot3}=\frac1{72}=T_3-T_2,\quad \frac{40-36}{40\cdot3\cdot2}=\frac1{60}=T_4-T_3.
		\]
		For $\Gamma_{126}$, the complete tripartite graph $K_{3,3,3}$ is walk-regular with nonzero Laplacian eigenvalues $6$ (six times) and $9$ (twice). The modules $M_1,\dots,M_5$ have $c_i=9,36,18,12,18$ and
		\begin{gather*}
			T_1=\frac{32}{33\cdot42}=\frac{16}{693},\qquad T_2=\frac{8}{9\cdot36}=\frac2{81},\qquad T_3=\frac{2}{3\cdot18}=\frac1{27},\\
			T_4=\frac19\Bigl(\frac{6}{18}+\frac{2}{21}\Bigr)=\frac1{21},\qquad T_5=\frac{8}{9\cdot18}=\frac4{81},
		\end{gather*}
		and for the edges $\{1,2\},\dots,\{4,5\}$ we have $A_{ij}=33,42,45,54$, so that
		\begin{gather*}
			\frac{60}{126\cdot33\cdot9}=\frac{10}{6237}=T_2-T_1,\qquad \frac{42}{126\cdot27}=\frac1{81}=T_3-T_2,\\
			\frac{36}{126\cdot27}=\frac2{189}=T_4-T_3,\qquad \frac{18}{126\cdot81}=\frac1{567}=T_5-T_4.
		\end{gather*}
		By Corollary~\ref{cor:tree}, both graphs are resistance-regular, and the stated values of $\Lp_{vv}$ follow from the formula in Corollary~\ref{cor:tree}. The degrees of the vertices in the first half of $\Gamma_{40}$ are $8+6=14$, $9+3=12$, $6+2=8$ and $3+2=5$; in $\Gamma_{126}$ they are $32+9=41$, $33+3=36$, $9+9=18$, $6+3+9=18$ and $9+9=18$.
	\end{proof}
	
	\subsection{Proof of Theorem~\ref{thm:far}}
	
	\begin{proof}[Proof of Theorem~\ref{thm:far}]
		(i) Let $t\ge2$ and $G=C_6(t,t^2)$. By Lemma~\ref{lem:c6} with $x=t$ and $y=t^2=x^2$, the graph $G$ is resistance-regular, and its common diagonal entry simplifies to
		\[
		\Lp_{vv}=\frac{16t^4+16t^3+2t+1}{8t^4(2t+1)^2}.
		\]
		The graph $G$ is connected and has $n=2x+4y=4t^2+2t$ vertices. The $2t$ vertices of the independent modules have degree $2y=2t^2$, and the $4t^2$ vertices of the clique modules have degree $x+2y-1=2t^2+t-1$. These degrees differ by $t-1\ge1$, and since $4n+1=(4t+1)^2$, we have $t-1=(\sqrt{4n+1}-5)/4$.
		
		(ii) By Lemma~\ref{lem:paths}, $\Gamma_{126}$ has exactly three distinct degrees and $\Gamma_{40}$ has exactly four.
		
		(iii) Let $G$ be a connected resistance-regular graph on $n\ge3$ vertices. Let $u$ be a vertex of degree $\Delta$, put $D=\Delta+1$, and let $F=V\setminus(\{u\}\cup N(u))$, so that $|F|=n-D$. By Theorem~\ref{thm:zwb}(iii),
		\[
		R=R_u=\Bigl(2-\frac2n\Bigr)+\sum_{w\in F}\rres{u}{w}.
		\]
		For $w\in F$ and $z\in\{u\}\cup N(u)$, the triangle inequality gives $\rres{w}{z}\ge\rres{u}{w}-\rres{u}{z}$. Summing over the $D$ vertices $z$ and using Theorem~\ref{thm:zwb}(iii) at $u$ once more, we get $R=R_w\ge D\rres{u}{w}-(2-2/n)$. Summing over $w\in F$ gives $(n-D)\bigl(R+2-2/n\bigr)\ge D\bigl(R-2+2/n\bigr)$, that is, $(2D-n)R\le 2(n-1)$. On the other hand, at a vertex $w_0$ of degree $\delta$, the third inequality of Lemma~\ref{lem:dirichlet} gives $\rres{w_0}{z}\ge\frac1{\delta+1}+\frac1{D}$ for every $z\ne w_0$, so $R\ge(n-1)\bigl(\frac1{\delta+1}+\frac1D\bigr)$. If $2D>n$, combining the two bounds gives
		\[
		\frac1{\delta+1}\le\frac{2}{2D-n}-\frac1D=\frac{n}{D(2D-n)},
		\]
		that is, $(\delta+1)\,n\ge(\Delta+1)(2\Delta+2-n)$. Then $\Delta-\delta=D-(\delta+1)\le D-D(2D-n)/n=2D(n-D)/n\le n/2$, with equality only if $D=n/2$, which is excluded. If $2D\le n$, then $\Delta-\delta\le\Delta-1<D\le n/2$. Hence $\Delta-\delta<n/2$ in all cases.
		
		Now suppose that $G$ has exactly two distinct degrees $\delta<\Delta$, and put $g=\Delta-\delta$. Because $G$ is connected and has two degree classes, some vertex $x$ of degree $\Delta$ has $b\ge1$ neighbours of degree $\delta$. By Theorem~\ref{thm:zwb}(iii) and the second inequality of Lemma~\ref{lem:dirichlet},
		\[
		2-\frac2n=\sum_{u\sim x}\rres{x}{u}\ge\sum_{u\sim x}\frac{\Delta+d_u-2}{\Delta d_u-1}.
		\]
		A term with $d_u=\Delta$ equals $2/(\Delta+1)$, and a term with $d_u=\delta$ equals $2/(\Delta+1)+g(\Delta-1)/\bigl((\Delta\delta-1)(\Delta+1)\bigr)$. Hence
		\[
		\frac{b\,g\,(\Delta-1)}{\Delta\delta-1}\le(\Delta+1)\Bigl(2-\frac2n\Bigr)-2\Delta=\frac{2(n-1-\Delta)}{n}.
		\]
		Since $(\Delta-1)(\delta+1)-(\Delta\delta-1)=\Delta-\delta>0$, we have $(\Delta-1)/(\Delta\delta-1)>1/(\delta+1)$, and with $b\ge1$ we obtain
		\begin{equation}\label{eq:twodegrees}
			\Delta-\delta<\frac{2(\delta+1)(n-1-\Delta)}{n}.
		\end{equation}
		Substituting $\Delta=\delta+g$ and rearranging gives $g<f(\delta+1)$ with $f(z)=2z(n-z)/(n+2z)$, and the maximum of $f$ over $z>0$ is attained at $z=(\sqrt3-1)n/2$ and equals $(2-\sqrt3)n$.
	\end{proof}
	
	The proof of part (iii) gives slightly more: if $2(\Delta+1)>n$, then $(\delta+1)\,n\ge(\Delta+1)(2\Delta+2-n)$, and with two degrees \eqref{eq:twodegrees} holds. For $\Delta=n-1$ the first inequality gives $\delta=n-1$, which recovers the fact that a resistance-regular graph with a vertex adjacent to all others is complete. For $t=1$, $C_6(1,1)$ is the $6$-cycle. The graphs $C_6(2,4)$ and $C_6(3,9)$ have $20$ and $42$ vertices; the first is, as far as we know, the second smallest non-regular resistance-regular graph, after $\Fp(4)$ (see Table~\ref{tab:new}). The graph $\Gamma_{40}$ has $\Delta-\delta=9$, that is, $0.225\,n$.
	
	\subsection{Algebraic connectivity and the largest difference of degrees}
	
	Large differences of degrees also force the graph to have a bottleneck. Recall that the algebraic connectivity $\mu_1$ of $G$ is its smallest nonzero Laplacian eigenvalue.
	
	\begin{proposition}\label{prop:bottleneck}
		Every connected resistance-regular graph on $n\ge3$ vertices satisfies
		\[
		(\Delta-\delta)(n-2)\,\mu_1\le\delta\,(n-\Delta).
		\]
	\end{proposition}
	
	\begin{proof}
		For a vertex $v$ let $\pi_v=\frac1{d_v}\sum_{u\sim v}e_u$ be the uniform distribution on its neighbours. Since $Le_v=d_v(e_v-\pi_v)$, we have $(e_v-\pi_v)^\top\Lp(e_v-\pi_v)=e_v^\top L\Lp Le_v/d_v^2=1/d_v$ and $e_v^\top\Lp(e_v-\pi_v)=e_v^\top\Lp Le_v/d_v=(1-1/n)/d_v$. Expanding $\pi_v^\top\Lp\pi_v$ with $\pi_v=e_v-(e_v-\pi_v)$ gives the identity
		\[
		\Lp_{vv}=\pi_v^\top\Lp\pi_v+\Bigl(1-\frac2n\Bigr)\frac1{d_v},
		\]
		valid in every connected graph. Let $x$ have degree $\Delta$ and $w$ degree $\delta$. If $G$ is resistance-regular, subtracting the identities at $x$ and $w$ gives $\pi_x^\top\Lp\pi_x-\pi_w^\top\Lp\pi_w=(1-2/n)(1/\delta-1/\Delta)$. Since $\Lp\one=0$ and $\pi_x-\one/n\perp\one$, we have $\pi_x^\top\Lp\pi_x=(\pi_x-\one/n)^\top\Lp(\pi_x-\one/n)\le\|\pi_x-\one/n\|^2/\mu_1=(1/\Delta-1/n)/\mu_1$, while $\pi_w^\top\Lp\pi_w\ge0$. Therefore $(1-2/n)(\Delta-\delta)/(\delta\Delta)\le(n-\Delta)/(n\Delta\mu_1)$, which rearranges to the claim.
	\end{proof}
	
	Write $g(n)$ for the largest value of $\Delta-\delta$ over all connected resistance-regular graphs on $n$ vertices. Theorem~\ref{thm:far} gives
	\[
	\frac{\sqrt{4n+1}-5}{4}\le g(n)<\frac n2\qquad\text{for }n=4t^2+2t,\ t\ge2,
	\]
	and by Proposition~\ref{prop:bottleneck}, resistance-regular graphs whose degrees differ substantially must have small algebraic connectivity. The order of growth of $g(n)$ is open; see Section~\ref{sec:consequences}.
	
	\subsection{Further examples}
	
	Additional examples with three distinct degrees can be obtained by blowing up other trees and block graphs, the smallest of which after $\Gamma_{126}$ has $192$ vertices (described in Appendix~\ref{app:examples}). Searching over blow-ups of all connected base graphs on at most seven vertices with walk-regular modules yielded no example with three or more degrees on fewer than $40$ vertices, and an exhaustive check of all $3.1\times10^{8}$ connected tri-circulants over $\Z_k$ for $k\le7$ (graphs admitting a three-orbit semiregular cyclic automorphism group) produced no non-regular resistance-regular graph at all. Whether the number of distinct vertex degrees in a resistance-regular graph can be arbitrarily large remains open.
	
	The constructions of Sections~\ref{sec:modules} and~\ref{sec:far}, together with searches over module blow-ups and over graphs with a two-orbit abelian group of automorphisms, produced many further examples. We collected $138$ pairwise non-isomorphic non-regular resistance-regular graphs on at most $300$ vertices, among them examples of odd order, each verified in exact rational arithmetic; they are listed in the ancillary files. Table~\ref{tab:new} shows a selection: the smallest examples we know with degree differences $2$, $3$ and $4$, the two smallest examples of odd order, and the graphs with three and four distinct degrees.
	
	\begin{table}[!htbp]
		\centering
		\caption{A selection of further non-regular resistance-regular graphs; degrees are listed with the number of vertices of each degree as an exponent. The graphs on $24$, $55$, $60$ and $192$ vertices are described in Appendix~\ref{app:examples}, the graph on $33$ vertices is given in the ancillary files, and all are available in \texttt{graph6} format there.}
		\label{tab:new}
		\begin{adjustbox}{max width=\textwidth}
			\begin{tabular}{@{}r r l l l @{\hspace{2em}} r r l l l@{}}
				\toprule
				$n$ & $|E|$ & degrees & $\Lp_{vv}$ & construction & $n$ & $|E|$ & degrees & $\Lp_{vv}$ & construction\\
				\midrule
				$20$ & $88$ & $8^4,\,9^{16}$ & $389/3200$ & $C_6(2,4)$ & $55$ & $260$ & $9^{30},\,10^{25}$ & $1162/9075$ & App.~\ref{app:examples}\\
				$24$ & $100$ & $7^8,\,9^{16}$ & $223/1440$ & App.~\ref{app:examples} & $60$ & $786$ & $23^{12},\,27^{48}$ & $7/150$ & App.~\ref{app:examples}\\
				$33$ & $187$ & $11^{22},\,12^{11}$ & $694/7623$ & ancillary files & $72$ & $1248$ & $32^8,\,35^{64}$ & $5129/165888$ & $C_6(4,16)$\\
				$40$ & $232$ & $5^4,\,8^6,\,12^{12},\,14^{18}$ & $99/400$ & $\Gamma_{40}$ & $126$ & $2055$ & $18^{42},\,36^{18},\,41^{66}$ & $1193/15876$ & $\Gamma_{126}$\\
				$42$ & $414$ & $18^6,\,20^{36}$ & $1735/31752$ & $C_6(3,9)$ & $192$ & $5412$ & $35^{12},\,56^{72},\,59^{108}$ & $1397/36864$ & App.~\ref{app:examples}\\
				\bottomrule
			\end{tabular}
		\end{adjustbox}
	\end{table}

	\section{Resistance rows and the degree}\label{sec:rows}
	
	Our resistance-regular graphs show that the sum $R_v$ of the multiset $\{\!\{\rres{v}{u}:u\in V\}\!\}$ does not determine $d_v$. Even the entire multiset $\{\!\{\rres{v}{u}:u\in V\}\!\}$ at a single vertex $v$ does not determine $d_v$, as happens already on ten vertices.
	
	\begin{proposition}\label{prop:tenvertices}
		Up to isomorphism, there are exactly two connected graphs on at most ten vertices that contain a pair of vertices $x,y$ of different degrees with $\{\!\{\rres{x}{u}:u\in V\}\!\}=\{\!\{\rres{y}{u}:u\in V\}\!\}$. Both have ten vertices:
		\begin{enumerate}[label=(\alph*),leftmargin=2.5em]
			\item $K_{2,3}$ with one vertex $x$ of the $3$-element part joined to a new vertex $y$ that is additionally adjacent to two pendant vertices and the two vertices of a triangle; here $d_x=3$, $d_y=5$, and both multisets of nonzero resistances equal $\{\!\{\tfrac23,\tfrac23,1,1,1,\tfrac53,\tfrac53,2,2\}\!\}$;
			\item $K_{2,3}$ with one vertex $x$ of the $2$-element part joined to the central vertex $y$ of a bowtie (two triangles sharing $y$); here $d_x=4$, $d_y=5$, and both multisets of nonzero resistances equal $\{\!\{\tfrac23,\tfrac23,\tfrac23,\tfrac23,1,\tfrac53,\tfrac53,\tfrac53,\tfrac53\}\!\}$.
		\end{enumerate}
		In both graphs, the second round of resistance-based colour refinement separates $x$ from $y$, and the stable colouring coincides with the automorphism orbit partition.
	\end{proposition}
	
	\begin{proof}
		We generated all $11\,989\,762$ connected graphs on $n\le10$ vertices using \texttt{geng}, computed their effective resistances exactly as in Section~\ref{sec:enum} (as the fractions $\tau\rres{u}{v}/\tau$ with integer numerators), and compared the sorted rows of $R$ within each graph; every candidate pair was verified again in exact rational arithmetic. For $n\le9$, every pair of vertices with identical resistance multisets has the same degree, whereas for $n=10$ only the two graphs in (a) and (b) contain a pair of vertices of unequal degree with the same resistance multiset.
	\end{proof}
	
	Although two individual vertices of different degrees can share the same resistance multiset, requiring \emph{every} vertex of $G$ to have the same resistance multiset, so that even the stable resistance-based Weisfeiler--Leman colouring cannot distinguish any two vertices, is much more restrictive.
	
	\begin{conjecture}\label{conj:rows}
		Let $G$ be a connected graph in which the multiset $\{\!\{\rres{v}{u}:u\in V\}\!\}$ is the same for every vertex $v\in V$. Then $G$ is regular.
	\end{conjecture}
	
	The hypothesis of Conjecture~\ref{conj:rows} means that all rows of the resistance matrix $R$ are permutations of one another. Any such graph is automatically resistance-regular, and in a resistance-regular graph the second resistance moment satisfies
	\begin{equation}\label{eq:secondmoment}
		\sum_{u\in V}\rres{v}{u}^2=n(\Lp_{vv})^2+2\Lp_{vv}\tr\Lp+\sum_{u\in V}(\Lp_{uu})^2+4\bigl((\Lp)^2\bigr)_{vv}.
	\end{equation}
	If all rows of $R$ are permutations of one another, the left-hand side of \eqref{eq:secondmoment} does not depend on $v$, and hence neither does $\bigl((\Lp)^2\bigr)_{vv}$. Consequently, by Theorem~\ref{thm:forced}(i), Conjecture~\ref{conj:rows} holds whenever $L$ has $m\le3$ distinct nonzero eigenvalues (and the second moment rules out all sixteen graphs of Table~\ref{tab:examples}; see Section~\ref{sec:moment}). Using the Dirichlet bounds of Lemma~\ref{lem:dirichlet} and the following maximum principle, we now prove that any connected graph whose resistance matrix has rows that are permutations of one another can have at most two consecutive degrees.
	
	\begin{lemma}[nearest vertices are neighbours]\label{lem:nearest}
		Let $G$ be a connected resistance-regular graph and let $v\in V$. Every vertex $u\ne v$ that minimises $\rres{v}{u}$ is adjacent to $v$.
	\end{lemma}
	
	\begin{proof}
		Because $\Lp_{uu}$ is constant across $u\in V$, formula~\eqref{eq:resistance} shows that $\rres{v}{u}$ is minimised over $u\ne v$ precisely when $\varphi(u)=\Lp_{uv}$ is maximised over $u\ne v$, where $\varphi=\Lp e_v$. For any vertex $w\ne v$, the Laplacian equation $(L\varphi)_w=-1/n$ can be written as
		\[
		\varphi(w)=\frac{1}{d_w}\sum_{u\sim w}\varphi(u)-\frac{1}{n d_w},
		\]
		so $\varphi(w)$ is strictly smaller than the maximum of $\varphi$ over the neighbours of $w$. If a vertex $u^*\in V\setminus\{v\}$ maximising $\varphi$ on $V\setminus\{v\}$ were not adjacent to $v$, then all neighbours of $u^*$ would lie in $V\setminus\{v\}$ and at least one of them would have a strictly larger value of $\varphi$, contradicting the choice of $u^*$.
	\end{proof}
	
	\begin{proposition}\label{prop:nearregular}
		Let $G$ be a connected graph on $n$ vertices in which all rows of the resistance matrix $R$ are permutations of one another, and let $\Delta$ and $\delta$ be the maximum and minimum degrees of $G$. Then $\Delta-\delta\le1$. Moreover, if $G$ is not regular, so that its vertex degrees are $k$ and $k+1$ for some $k\ge2$, then $n\ge2k+5$ when $k\ge3$, and $n\ge10$ when $k=2$.
	\end{proposition}
	
	\begin{proof}
		Because the rows of $R$ are permutations of one another, $G$ is resistance-regular and the minimum positive entry $\rho=\min_{u\ne v}\rres{v}{u}$ in row $v$ does not depend on $v$. Choose a vertex $v$ of maximum degree $\Delta$. By Theorem~\ref{thm:zwb}(iii), the $\Delta$ edge resistances incident to $v$ sum to $2-2/n$, so their minimum satisfies
		\[
		\rho\le\frac{2-2/n}{\Delta}<\frac{2}{\Delta}.
		\]
		Now let $w\in V$ be any vertex, and choose $u\ne w$ with $\rres{w}{u}=\rho$. By the third inequality of Lemma~\ref{lem:dirichlet}, $\rho\ge\frac{1}{d_w+1}+\frac{1}{d_u+1}\ge\frac{1}{d_w+1}+\frac{1}{\Delta+1}$. Combining the upper and lower bounds on $\rho$ gives
		\[
		\frac{1}{d_w+1}<\frac{2}{\Delta}-\frac{1}{\Delta+1}=\frac{\Delta+2}{\Delta(\Delta+1)},\qquad\text{so}\qquad
		d_w+1>\frac{\Delta(\Delta+1)}{\Delta+2}=\Delta-1+\frac{2}{\Delta+2}.
		\]
		Since $d_w$ is an integer, $d_w\ge\Delta-1$, proving that $\Delta-\delta\le1$.
		
		Suppose now that $G$ is not regular, with degrees $k$ and $k+1$; by Proposition~\ref{prop:necessary}, $k\ge2$. Let $w$ be a vertex of degree $k$ and let $u\ne w$ achieve $\rres{w}{u}=\rho$. By Lemma~\ref{lem:nearest}, $u$ is adjacent to $w$, and since $d_u\in\{k,k+1\}$, the adjacent bound of Lemma~\ref{lem:dirichlet} gives
		\[
		\rho\ge\min\Bigl\{\frac{2k-2}{k^2-1},\,\frac{2k-1}{k^2+k-1}\Bigr\}=\frac{2k-1}{k^2+k-1}.
		\]
		Comparing this with the upper bound $\rho\le(2-2/n)/(k+1)$ obtained at a vertex of degree $\Delta=k+1$ gives $(2k-1)(k+1)\le(2-2/n)(k^2+k-1)$, which rearranges to
		\[
		n\ge\frac{2(k^2+k-1)}{k-1}=2k+4+\frac{2}{k-1}.
		\]
		Since $n$ is an integer, this yields $n\ge2k+5$ for $k\ge3$ and $n\ge10$ for $k=2$.
	\end{proof}
	
	Proposition~\ref{prop:nearregular} should be compared with Theorem~\ref{thm:far}. When only the row sums of $R$ are required to coincide, the degrees can differ by $(\sqrt{4n+1}-5)/4$ and there can be three or four distinct degrees; when the full rows are required to coincide as multisets, the degrees differ by at most one.
	
	To rule out the remaining two-degree case $\{k,k+1\}$ of Conjecture~\ref{conj:rows} when $m\ge4$, spectral moments alone do not suffice because the constancy of $\Lp_{vv}$ and $\bigl((\Lp)^2\bigr)_{vv}$ only constrains the moments of orders $0,-1,-2$ of the local spectral measure. Any full proof must therefore use the $0/{-1}$ off-diagonal structure of the Laplacian; Example~\ref{ex:weighted} below shows that a related statement fails as soon as the edges carry weights. One local identity that brings the neighbourhood $N(v)$ directly into play expresses the sum of the effective resistances \emph{between} the neighbours of $v$ as a linear function of $d_v$.
	
	\begin{lemma}[resistances among the neighbours of a vertex]\label{lem:neighbours}
		For every vertex $v$ of a connected graph $G$, with $\sigma_v=\sum_{u\sim v}\rres{u}{v}$,
		\[
		\sum_{u\in N(v)}\sum_{w\in N(v)}\rres{u}{w}=2d_v(\sigma_v-1).
		\]
		In particular, in a resistance-regular graph on $n$ vertices, the resistances between the neighbours of $v$ sum to $2d_v(1-2/n)$.
	\end{lemma}
	
	\begin{proof}
		Write $a_u=\Lp_{uu}$. The $(v,w)$ entry of $L\Lp=I-J/n$ reads $d_v\Lp_{vw}-\sum_{u\sim v}\Lp_{uw}=\delta_{vw}-1/n$; taking $w=v$ gives $\sum_{u\sim v}\Lp_{uv}=d_va_v-1+1/n$, and summing the general case over $w\sim v$ gives $\sum_{u,w\in N(v)}\Lp_{uw}=d_v\sum_{w\sim v}\Lp_{vw}+d_v/n=d_v(d_va_v-1+1/n)+d_v/n$. By \eqref{eq:resistance}, $\sum_{u,w\in N(v)}\rres{u}{w}=2d_v\sum_{u\sim v}a_u-2\sum_{u,w\in N(v)}\Lp_{uw}=2d_v\bigl(\sum_{u\sim v}a_u-d_va_v+1-2/n\bigr)$, while $\sigma_v=\sum_{u\sim v}(a_u+a_v-2\Lp_{uv})=\sum_{u\sim v}a_u-d_va_v+2-2/n$. Comparing the two expressions gives the identity, and the last statement follows from Theorem~\ref{thm:zwb}(iii).
	\end{proof}
	
	Computationally, Conjecture~\ref{conj:rows} holds for all connected graphs on at most $10$ vertices: for all $11\,989\,762$ of them, the stable resistance-based colouring coincides with the partition into automorphism orbits (an exact statement, since the resistances were computed exactly, the orbits refine the stable colouring and hash collisions can only merge colour classes), so a graph with a single colour class is vertex-transitive. It also holds for all $3.6\times10^{11}$ bi-Cayley graphs of Section~\ref{sec:search}. The hypothesis of the conjecture differs from walk-regularity even for regular graphs: the non-vertex-transitive walk-regular graphs on $12$ vertices, the Folkman graph and the Gray graph have distinct row multisets, and among about $33\,000$ small regular graphs (all cubic graphs on at most $16$ vertices, the $4$-regular graphs on at most $13$, the $5$-regular graphs on $10$ and $12$, and the $6$-regular graphs on $11$ and $12$ vertices) every graph whose resistance rows coincide as multisets is vertex-transitive.
	
	\subsection{The second resistance moment}\label{sec:moment}
	
	The proof of Proposition~\ref{prop:nearregular} uses only the smallest entry of each row of $R$. The examples suggest that a single further statistic of a row, its second moment, already determines the degree in every resistance-regular graph. For a vertex $v$ write $S_2(v)=\sum_{u\in V}\rres{v}{u}^2$. By \eqref{eq:secondmoment}, in a resistance-regular graph $S_2(v)$ differs from $4\bigl((\Lp)^2\bigr)_{vv}$ by a quantity that does not depend on $v$.
	
	\begin{conjecture}\label{conj:moment}
		Let $G$ be a connected resistance-regular graph and let $u,v\in V$. If $d_u<d_v$, then $S_2(u)<S_2(v)$; equivalently, $\bigl((\Lp)^2\bigr)_{uu}<\bigl((\Lp)^2\bigr)_{vv}$.
	\end{conjecture}
	
	If all rows of $R$ are permutations of one another, then $G$ is resistance-regular and $S_2$ is constant, so Conjecture~\ref{conj:moment} implies Conjecture~\ref{conj:rows}. Similarly, if both $\Lp$ and $(\Lp)^2$ have constant diagonal, then $G$ is resistance-regular and, by \eqref{eq:secondmoment}, $S_2$ is constant. Hence Conjecture~\ref{conj:moment} implies that a connected graph in which $\Lp$ and $(\Lp)^2$ both have constant diagonal is regular; in other words, the first two of the conditions in Theorem~\ref{thm:forced}(i) would force regularity, whatever the value of $m$. For $m\le2$ the first condition suffices (Theorem~\ref{thm:forced}(ii)), and for $m=3$ the two conditions are all of Theorem~\ref{thm:forced}(i). Since $S_2(v)$ is determined by the multiset $\{\!\{\rres{v}{u}:u\in V\}\!\}$, Conjecture~\ref{conj:moment} would also imply that in a resistance-regular graph the first round of the resistance-based Weisfeiler--Leman test already separates vertices of different degrees.
	
	The conjecture holds in every example known to us. We verified it in exact rational arithmetic for all graphs of Tables~\ref{tab:examples} and~\ref{tab:new}, for all $138$ graphs of the ancillary catalogue (for the graphs with three and four degrees, the inequality holds between every pair of degree classes), and for the $3654$ members of the four-parameter family of Appendix~\ref{app:fourparam} that our searches produced and that are listed in the ancillary files. It can also be proved in several cases.
	
	\begin{proposition}\label{prop:momentcases}
		Conjecture~\ref{conj:moment} holds in each of the following cases:
		\begin{enumerate}[label=(\alph*),leftmargin=2.5em]
			\item $L$ has exactly three distinct nonzero eigenvalues;
			\item $L$ has exactly four distinct nonzero eigenvalues, with sum $\sigma_1$, and $d_u+d_v+1<\sigma_1$ whenever $d_u<d_v$;
			\item $G=\Fp(p)$ for an even $p\ge4$;
			\item $G=C_6(t,t^2)$ for an integer $t\ge2$;
			\item $G$ is obtained by the construction of Proposition~\ref{prop:gap1} from a graph $H$ that satisfies Conjecture~\ref{conj:moment}.
		\end{enumerate}
	\end{proposition}
	
	\begin{proof}
		Case (a) is Corollary~\ref{cor:m3}. For (b), let $\sigma_1,\dots,\sigma_4$ be the elementary symmetric functions of the four distinct nonzero eigenvalues. By Lemma~\ref{lem:minpoly} with $m=4$ and $j=2$, $L^2-\sigma_1L+\sigma_2P-\sigma_3\Lp+\sigma_4(\Lp)^2=0$, where $P=I-J/n$. On the diagonal, with $\Lp_{vv}$ constant, this yields
		\[
		\sigma_4\Bigl(\bigl((\Lp)^2\bigr)_{vv}-\bigl((\Lp)^2\bigr)_{uu}\Bigr)=(d_v-d_u)(\sigma_1-1-d_u-d_v),
		\]
		which has the sign of $d_v-d_u$ under the stated condition. For (c), the diagonal of $(\Lp)^2$ at $(h,i)$ is the average of the $(i,i)$ entries of $(\widehat L_\chi^{\dagger})^2$ over all characters $\chi$, as in the proof of Lemma~\ref{lem:bicayley}. Evaluating these averages with the table of Lemma~\ref{lem:Fbicayley} shows that between a vertex $v\in Q_b$ (degree $2p-1$) and a vertex $w\in X_b$ (degree $2p-2$),
		\[
		\bigl((\Lp)^2\bigr)_{vv}-\bigl((\Lp)^2\bigr)_{ww}=\frac{7p-6}{8p^2(p-1)(p-2)}>0;
		\]
		for $\Fp(4)$, this gives $S_2(v)=371/128$ at vertices of degree $7$ and $1069/384$ at vertices of degree $6$. For (d), Proposition~\ref{prop:blowup} together with the values $y^{(0)}$ and $y^{(1)}$ in the proof of Lemma~\ref{lem:c6} gives, for $w$ in a clique module and $v$ in an independent module of $C_6(t,t^2)$,
		\[
		\bigl((\Lp)^2\bigr)_{ww}-\bigl((\Lp)^2\bigr)_{vv}=\frac{(t-1)(4t^2+2t+1)}{8t^6(2t+1)}>0.
		\]
		Case (e) follows from the second formula of Proposition~\ref{prop:gap1}, since the substitution preserves the order of the degrees.
	\end{proof}
	
	For graphs~14 and~15 of Table~\ref{tab:examples}, which have $m=4$, the nonzero eigenvalues $12$, $16$ and $12\pm4\sqrt5$ have sum $52>14+15+1$, in accordance with (b).
	
	The conjecture depends on the edges having equal weights. For a graph with positive edge weights, the weighted Laplacian, its pseudoinverse, the effective resistances and the weighted degrees are defined in the usual way, and Theorem~\ref{thm:zwb}(ii) and \eqref{eq:secondmoment} remain valid.
	
	\begin{example}\label{ex:weighted}
		Let $H$ be the house graph on the vertices $0,\dots,4$, consisting of the $4$-cycle $0,3,1,4$ and the roof vertex $2$ adjacent to $0$ and $4$, with edge weights $w_{02}=w_{24}=360$, $w_{03}=w_{14}=273$, $w_{04}=15$ and $w_{13}=455$. The weighted degrees of the vertices $0,1,2,3,4$ are $648$, $728$, $720$, $728$ and $648$, and every diagonal entry of the pseudoinverse of the weighted Laplacian equals $121/102375$, so $H$ is resistance-regular in the weighted sense. Nevertheless $S_2(2)-S_2(1)=79/1117935000>0$, although vertex $2$ has smaller weighted degree than vertex $1$.
	\end{example}
	
	Thus no argument that applies verbatim to weighted graphs can prove Conjecture~\ref{conj:moment}. The conjecture also has a probabilistic reading. Let $\eta_v$ be the hitting time of $v$ for the continuous-time random walk on $G$ that crosses each edge at rate $1$, started from a uniformly random vertex. The expected hitting time from $u$, $\psi(u)=n(\Lp_{vv}-\Lp_{uv})$, satisfies $(L\psi)(u)=1$ for $u\ne v$ and $\psi(v)=0$, and a similar computation for the second moment gives
	\[
	\mathbb E\,\eta_v=n\Lp_{vv},\qquad \operatorname{Var}\eta_v=2n\bigl((\Lp)^2\bigr)_{vv}+n^2(\Lp_{vv})^2.
	\]
	In a resistance-regular graph all vertices are thus equally hard to hit on average, and Conjecture~\ref{conj:moment} states that the hitting time of a vertex of larger degree has the larger variance. We intend to study Conjecture~\ref{conj:moment} further in a separate paper.

	\section{Consequences and open problems}\label{sec:consequences}
	
	Resistance-regularity appears under several names in combinatorics, network analysis and machine learning, and our results have the following consequences.
	\begin{itemize}[leftmargin=2em,itemsep=2pt]
		\item \emph{Resistance curvature.} The vertex curvature $p_v=1-\tfrac12\sum_{u\sim v}\rres{u}{v}$ of Devriendt and Lambiotte \cite{DevriendtLambiotte2022} is constant if and only if $G$ is resistance-regular. Theorem~\ref{thm:existence} therefore gives the first infinite family of non-regular graphs with constant resistance curvature, answering the question raised in \cite{DevriendtOttoliniSteinerberger2024}, where all known examples were vertex-transitive \cite{DevriendtOttoliniSteinerberger2024, SunYangXu2025}. By Theorem~\ref{thm:far}, the degrees of such a graph can differ by about $\sqrt n/2$, and there can be three or four of them.
		\item \emph{Current-flow closeness centrality.} The current-flow closeness of Brandes and Fleischer \cite{BrandesFleischer2005}, equivalent up to scaling to the information centrality of Stephenson and Zelen \cite{StephensonZelen1989}, is $(n-1)/R_v$. In a resistance-regular graph all vertices have the same current-flow closeness, even when their degrees differ, by $t-1$ in $C_6(t,t^2)$; by Proposition~\ref{prop:bottleneck}, a large difference of degrees then forces a small algebraic connectivity.
		\item \emph{Uniform spanning trees.} By \eqref{eq:kirchhoff} and Theorem~\ref{thm:zwb}(iii), in a uniformly random spanning tree of a resistance-regular graph the degree of every vertex has expectation $2(n-1)/n$. Higher moments still distinguish the degree classes: in $\Fp(4)$, the variance of the tree degree is about $0.682$ at a vertex of degree $6$ and $0.742$ at a vertex of degree $7$.
		\item \emph{Distance-based Weisfeiler--Leman tests.} Zhang et al.\ \cite{Zhang2023biconnectivity} studied the generalised distance Weisfeiler--Leman test, which refines the colour of $v$ by the multiset $\{\!\{(d(v,u),\operatorname{colour}(u)):u\in V\}\!\}$ for a graph metric $d$. For shortest-path distance the first round refines the degree partition; for effective resistance it is open whether even the stable colouring does. Flora et al.\ \cite{Flora2026} conjecture that the stable resistance-based colouring can fail to distinguish graphs that the classical test distinguishes; since the classical test refines the degree partition, a negative answer to our question would point towards their conjecture. For resistance-regular graphs, Conjecture~\ref{conj:moment} would imply that already the first round separates vertices of different degrees.
	\end{itemize}
	
	\paragraph{Open questions.}
	We close by recording the main questions left open by our results:
	\begin{enumerate}[leftmargin=2em,itemsep=2pt]
		\item \emph{Minimum order.} What is the smallest number of vertices of a non-regular resistance-regular graph? By Proposition~\ref{prop:enum} and Theorem~\ref{thm:existence}, it lies between $13$ and $16$.
		\item \emph{Three nonzero Laplacian eigenvalues.} Theorem~\ref{thm:forced}(ii) rules out $m\le2$, and Proposition~\ref{prop:sharp} provides examples with $m=4$. Does there exist a connected non-regular resistance-regular graph with $m=3$? By Theorem~\ref{thm:forced}(ii), it would have exactly two degrees, adding up to one less than the sum of the three eigenvalues.
		\item \emph{Difference of degrees.} What is the order of growth of the largest difference $g(n)$ between the maximum and the minimum degree of a resistance-regular graph on $n$ vertices? We know that $g(n)\ge(\sqrt{4n+1}-5)/4$ for $n=4t^2+2t$ and that $g(n)<n/2$, and $(2-\sqrt3)n$ is an upper bound when there are only two distinct degrees.
		\item \emph{Number of degrees and symmetry.} Is the number of distinct degrees of a resistance-regular graph bounded? Is there a non-regular resistance-regular graph with a trivial automorphism group? All examples known to us are built from modules or have a large automorphism group.
		\item \emph{Rows of the resistance matrix and diagonal powers.} Do Conjectures~\ref{conj:rows} and~\ref{conj:moment} hold? In particular, is every connected graph in which both $\Lp$ and $(\Lp)^2$ have constant diagonal regular? More generally, which combinations of resistance-based vertex invariants suffice to determine the degree of a vertex?
	\end{enumerate}
	
	\medskip
	\noindent\textbf{Declaration of generative AI and AI-assisted technologies.}
	During the preparation of this work, the authors used Claude Fable~5.1 to assist with writing the computational search scripts and brainstorming initial lines of inquiry, and Gemini Argon to assist with structuring, drafting and editing the manuscript text. After using these tools, the authors independently verified all mathematical proofs and exact computations, reviewed and edited the text as needed, and take full responsibility for the content of the publication.
	
	\appendix
	\section{The four-parameter family}\label{app:fourparam}
	
	Let $s,a,k\ge1$ and $h\ge2$ be integers. The graph $\mathcal G(s,a,k,h)$ has vertex set $\bigcup_{b<s}(Q_b\cup X_b)$, where the $Q_b$ are disjoint sets of size $a$ and $X_b=\{x_{b,c,t}:c<k,\ t<h\}$. Every $Q_b$ is a clique, every vertex of $Q_b$ is adjacent to every vertex of $X_b$, and $x_{b,c,t}\sim x_{b',c',t'}$ if and only if $c=c'$ and $t\ne t'$. Thus, for each $c<k$, the $sh$ vertices $x_{b,c,t}$ induce the complete multipartite graph whose $h$ parts $\{x_{b,c,t}:b<s\}$ have size $s$. The graph is connected and has $n=s(a+kh)$ vertices; a vertex of $Q_b$ has degree $a-1+kh$ and a vertex of $X_b$ has degree $a+s(h-1)$, so the two degrees differ by $kh-1-s(h-1)$. Comparing with Definition~\ref{def:family} shows that $\mathcal G(2,p,2,p/2)=\Fp(p)$. In this appendix only, we write
	\[
	E=a+kh,\qquad F=a+s(h-1),\qquad H=a+sh,\qquad D=khs(h-1).
	\]
	
	\begin{proposition}\label{prop:Gspectrum}
		The Laplacian eigenvalues of $\mathcal G(s,a,k,h)$ are $0$; $E$ with multiplicity $s(a-1)+1$; $a$ with multiplicity $k-1$; $H$ with multiplicity $k(h-1)$; $F$ with multiplicity $(s-1)(kh-1)$; and the two roots $\mu_\pm$ of $\mu^2-(kh+F)\mu+D$, each with multiplicity $s-1$.
	\end{proposition}
	
	\begin{proof}
		We decompose $\R^V$ into mutually orthogonal subspaces invariant under $L$:
		\begin{enumerate}[label=(\arabic*),leftmargin=2.5em]
			\item vectors supported on a single $Q_b$ whose entries sum to zero (dimension $s(a-1)$);
			\item vectors supported on $\bigcup_bX_b$ whose value at $x_{b,c,t}$ depends only on $c$ and sums to zero over $c$ (dimension $k-1$);
			\item vectors supported on $\bigcup_bX_b$ whose value at $x_{b,c,t}$ depends only on $(c,t)$ and sums to zero over $t$ for every $c$ (dimension $k(h-1)$);
			\item vectors supported on $\bigcup_bX_b$ whose values sum to zero over $b$ for every $(c,t)$ and over $(c,t)$ for every $b$ (dimension $(s-1)(kh-1)$);
			\item vectors that are constant on every $Q_b$ and on every $X_b$ (dimension $2s$).
		\end{enumerate}
		On (1), $L$ acts as $E$ by Lemma~\ref{lem:module}, since $Q_b$ is a clique module with $kh$ outside neighbours. At a vertex $q\in Q_b$, a vector $f$ in (2), (3) or (4) gives $(Lf)(q)=-\sum_{x\in X_b}f(x)=0$, by the zero-sum conditions. At a vertex $x=x_{b,c,t}$, $(Lf)(x)=(a+s(h-1))f(x)-\sum f(x_{b',c,t'})$, the sum running over $b'<s$ and $t'\ne t$; this sum equals $s(h-1)f(x)$ in (2), $-sf(x)$ in (3) and $0$ in (4), so $L$ acts as $a$, $H$ and $F$ respectively. Finally, if $f\in(5)$ takes the value $\alpha_b$ on $Q_b$ and $\beta_b$ on $X_b$, then $Lf$ takes the values $kh(\alpha_b-\beta_b)$ on $Q_b$ and $a(\beta_b-\alpha_b)+s(h-1)\beta_b-(h-1)\sum_{b'}\beta_{b'}$ on $X_b$. On the vectors with $\alpha_b,\beta_b$ independent of $b$ this is the matrix $\bigl(\begin{smallmatrix}kh&-kh\\-a&a\end{smallmatrix}\bigr)$, with eigenvalues $0$ and $E$, and on the vectors with $\sum_b\alpha_b=\sum_b\beta_b=0$ it acts on every pair $(\alpha_b,\beta_b)$ by the matrix
		\[
		M=\begin{pmatrix}kh&-kh\\-a&F\end{pmatrix},\qquad \det M=kh(F-a)=D,
		\]
		whose eigenvalues are $\mu_\pm$, each with multiplicity $s-1$. The dimensions add up to $s(a-1)+(k-1)+k(h-1)+(s-1)(kh-1)+2s=n$.
	\end{proof}
	
	\begin{proposition}\label{prop:Gdiag}
		For $q\in Q_b$ and $x\in X_b$,
		\begin{align*}
			\Lp_{qq}&=\frac{1}{sa}\Bigl(\frac{s(a-1)}{E}+\frac{kh}{E^2}+\frac{(s-1)F}{D}\Bigr),\\
			\Lp_{xx}&=\frac{1}{skh}\Bigl(\frac{k-1}{a}+\frac{k(h-1)}{H}+\frac{(s-1)(kh-1)}{F}+\frac{a}{E^2}+\frac{s-1}{s(h-1)}\Bigr),
		\end{align*}
		and
		\[
		\Lp_{qq}-\Lp_{xx}=\frac{N(s,a,k,h)}{ahks\,(a+hk)(a+hs)(a+hs-s)},
		\]
		where
		\begin{align*}
			N(s,a,k,h)={}&2(s-1)a^3+s\bigl(h^2ks-h^2k^2-hks+hk+3hs-3h-s+1\bigr)a^2\\
			&+hs\bigl(h^2ks^2-h^2k^2s-hks^2+hk+hs^2-hs+ks-k-s^2+s\bigr)a\\
			&-h^2ks^2(h-1)(k-1).
		\end{align*}
		Consequently, $\mathcal G(s,a,k,h)$ is resistance-regular if and only if $N(s,a,k,h)=0$.
	\end{proposition}
	
	\begin{proof}
		Since the subspaces in the proof of Proposition~\ref{prop:Gspectrum} are orthogonal and invariant, $\Lp_{vv}$ is the sum over them of $e_v^\top\Lp z_v$, where $z_v$ is the orthogonal projection of $e_v$ onto the subspace; on a subspace on which $L$ acts as a scalar $\mu$, this is $\|z_v\|^2/\mu$. For $q$, the squared norms of the projections onto (1)--(4) are $1-1/a$, $0$, $0$, $0$. For $x$ they are $0$, $(k-1)/(skh)$, $(h-1)/(sh)=k(h-1)/(skh)$ and $(s-1)(kh-1)/(skh)$, the last because the vectors in (4) are those that sum to zero over $b$, minus those that are constant on every $X_b$. In (5), the eigenvector for $E$ orthogonal to $\one$ is the vector $u$ with value $kh$ on every $Q_b$ and $-a$ on every $X_b$; since $\|u\|^2=sakhE$, it contributes $kh/(saE^2)$ at $q$ and $a/(skhE^2)$ at $x$. The projection of $e_q$ onto the remaining part of (5) has $\alpha_{b'}=(\delta_{bb'}-1/s)/a$ and $\beta=0$, and applying $\Lp$ means applying $M^{-1}=D^{-1}\bigl(\begin{smallmatrix}F&kh\\a&kh\end{smallmatrix}\bigr)$ to every pair; the value at $q$ is therefore $(1-1/s)(M^{-1})_{11}/a=(s-1)F/(saD)$. Similarly, the projection of $e_x$ has $\alpha=0$ and $\beta_{b'}=(\delta_{bb'}-1/s)/(kh)$, and it contributes $(1-1/s)(M^{-1})_{22}/(kh)=(s-1)/(sD)$. Adding up gives the two formulas, and bringing their difference to a common denominator, a routine computation that the ancillary files carry out symbolically, gives the last one. The denominator is positive, so $\Lp$ has constant diagonal if and only if $N(s,a,k,h)=0$.
	\end{proof}
	
	The same decomposition gives closed formulas for the diagonal of $(\Lp)^2$; they are checked in the ancillary scripts, and we used them to check Conjecture~\ref{conj:moment} for the members of the family listed below.
	
	As a polynomial in $k$, $N$ has negative leading coefficient $-h^2s\bigl(a^2+ahs+s(h-1)\bigr)$ and positive constant term $a(s-1)(a+hs)(2a+hs-s)$, so for given $s\ge2$, $a$ and $h$ there is exactly one positive real $k$ with $N=0$. Integer solutions are rarer, but they include a two-parameter family.
	
	\begin{corollary}\label{cor:Tfam}
		For integers $h\ge2$ and $m\ge0$, put $s=hm+2$, $k=(h+1)m+2$ and $a=sh(k-1)$. Then $\mathcal G(s,a,k,h)$ is a connected resistance-regular graph on $n=hs\bigl(h(h+1)m^2+(4h+3)m+4\bigr)$ vertices whose two degrees differ by $2hm+1=2s-3$. For $m=0$ it is $\Fp(2h)$.
	\end{corollary}
	
	\begin{proof}
		A direct expansion gives the polynomial identity $N\bigl(s,sh(k-1),k,h\bigr)=h^2ks^3(k-1)(hk-1)\bigl(s(h+1)-hk-2\bigr)$, and $s(h+1)-hk-2=0$ for the given $s$ and $k$. The remaining statements follow from the formulas for $n$ and for the degrees.
	\end{proof}
	
	Since $n$ is of order $m^3$ for fixed $h$, the degrees in Corollary~\ref{cor:Tfam} differ by an amount of order $n^{1/3}$, less than in Theorem~\ref{thm:far}(i). Table~\ref{tab:Gfam} lists, for each difference of degrees, the smallest member of the family found in a computer search over all parameters with $s\le4000$ and difference at most $12$ in absolute value; a negative difference means that the vertices of $X_b$ have the larger degree. In this range, $2$ is the only even difference that occurs. The family is infinite, since it contains all $\Fp(p)$ and all graphs of Corollary~\ref{cor:Tfam}; our searches, which also covered other ranges of parameters, produced a list of $3654$ members, given in the ancillary files, and each of them satisfies Conjecture~\ref{conj:moment}.
	
	\begin{table}[!htbp]
		\centering
		\caption{The smallest members of $\mathcal G(s,a,k,h)$ found for each difference of degrees $kh-1-s(h-1)$; the member with difference $1$ is $\Fp(4)$.}
		\label{tab:Gfam}
		\begin{adjustbox}{max width=\textwidth}
			\begin{tabular}{@{}r l r l @{\hspace{2.5em}} r l r l@{}}
				\toprule
				difference & $(s,a,k,h)$ & $n$ & degrees & difference & $(s,a,k,h)$ & $n$ & degrees\\
				\midrule
				$-1$ & $(6,9,5,6)$ & $234$ & $38,39$ & $5$ & $(4,32,5,2)$ & $168$ & $36,41$\\
				$1$ & $(2,4,2,2)$ & $16$ & $6,7$ & $7$ & $(5,75,6,3)$ & $465$ & $85,92$\\
				$2$ & $(33,33,18,2)$ & $2277$ & $66,68$ & $9$ & $(6,84,8,2)$ & $600$ & $90,99$\\
				$3$ & $(4,72,4,10)$ & $448$ & $108,111$ & $11$ & $(7,245,8,5)$ & $1995$ & $273,284$\\
				\bottomrule
			\end{tabular}
		\end{adjustbox}
	\end{table}
	
	\section{Connection sets and ancillary files}\label{app:examples}
	
	Each of the sixteen non-regular resistance-regular graphs in Table~\ref{tab:examples} is realised below as $\mathrm{BiCay}(H;S_{11},S_{22},S_{12})$ for a group $H$. Elements of $\Z_a\times\Z_b$ are written as pairs, elements of $\Z_2\times\Z_2\times\Z_4$ as triples, elements of the dihedral group $D_k=\langle r,s\mid r^k=s^2=1,\ srs=r^{-1}\rangle$ as $r^i$ or $r^is$, and elements of the dicyclic group $\mathrm{Dic}_3=\langle a,x\mid a^6=1,\ x^2=a^3,\ xax^{-1}=a^{-1}\rangle$ as $a^i$ or $a^ix$. For non-abelian $H$, edges are given by right multiplication: $(h,1)\sim(hs,1)$ for $s\in S_{11}$, $(h,2)\sim(hs,2)$ for $s\in S_{22}$, and $(h,1)\sim(hs,2)$ for $s\in S_{12}$.
	
	\begin{enumerate}[leftmargin=2em,itemsep=2pt]\sloppy
		\item $\Z_2\times\Z_4$: $S_{11}=\{(0,2),\allowbreak(1,0)\}$, \allowbreak $S_{22}=\{(0,1),\allowbreak(0,2),\allowbreak(0,3)\}$, \allowbreak $S_{12}=\{0\}\times\Z_4$.
		\item $\mathrm{Dic}_3$: $S_{11}=\{x,a^3x\}$, \allowbreak $S_{22}=\{a^3,ax,a^4x\}$, \allowbreak $S_{12}=\{1,a^3,ax,a^4x\}$.
		\item $\Z_{12}$: $S_{11}=\{\pm1,\pm2\}$, \allowbreak $S_{22}=\{2,4,6,8,10\}$, \allowbreak $S_{12}=\{0,6\}$.
		\item $\Z_2\times\Z_6$: $S_{11}=\{(0,1),\allowbreak(0,2),\allowbreak(0,4),\allowbreak(0,5)\}$, \allowbreak $S_{22}=\{(0,2),\allowbreak(0,4),\allowbreak(1,0),\allowbreak(1,2),\allowbreak(1,4)\}$, \allowbreak $S_{12}=\{(0,0),\allowbreak(0,2),\allowbreak(0,4),\allowbreak(1,0),\allowbreak(1,2),\allowbreak(1,4)\}$.
		\item $\Z_2\times\Z_8$: $S_{11}=\{(0,1),\allowbreak(0,4),\allowbreak(0,7),\allowbreak(1,1),\allowbreak(1,7)\}$, \allowbreak $S_{22}=\{(0,2),\allowbreak(0,6),\allowbreak(1,0),\allowbreak(1,2),\allowbreak(1,4),\allowbreak(1,6)\}$, \allowbreak $S_{12}=\{(0,0),\allowbreak(0,4)\}$.
		\item $\Z_4\times\Z_4$: $S_{11}=\{(0,1),\allowbreak(0,3),\allowbreak(1,0),\allowbreak(2,2),\allowbreak(3,0)\}$, \allowbreak $S_{22}=\{(0,1),\allowbreak(0,2),\allowbreak(0,3),\allowbreak(2,0),\allowbreak(2,1),\allowbreak(2,3)\}$, \allowbreak $S_{12}=\{(0,0),\allowbreak(0,2),\allowbreak(2,0),\allowbreak(2,2)\}$.
		\item $\Z_4\times\Z_4$: $S_{11}=\{(0,1),\allowbreak(0,3),\allowbreak(1,0),\allowbreak(1,1),\allowbreak(3,0),\allowbreak(3,3)\}$, \allowbreak $S_{22}=\{(0,1),\allowbreak(0,2),\allowbreak(0,3),\allowbreak(2,0),\allowbreak(2,1),\allowbreak(2,2),\allowbreak(2,3)\}$, \allowbreak $S_{12}=\{(0,0),\allowbreak(0,2),\allowbreak(2,1),\allowbreak(2,3)\}$.
		\item $\Z_4\times\Z_4$: $S_{11}$ and $S_{22}$ as in graph~7, \allowbreak $S_{12}=\{(0,0),\allowbreak(0,2),\allowbreak(2,0),\allowbreak(2,2)\}$.
		\item $\Z_4\times\Z_4$: $S_{11}=\{(0,1),\allowbreak(0,2),\allowbreak(0,3),\allowbreak(1,0),\allowbreak(2,0),\allowbreak(3,0)\}$, \allowbreak $S_{22}=\{(0,2),\allowbreak(1,1),\allowbreak(1,3),\allowbreak(2,0),\allowbreak(2,2),\allowbreak(3,1),\allowbreak(3,3)\}$, \allowbreak $S_{12}=\{(0,0),\allowbreak(1,1),\allowbreak(2,2),\allowbreak(3,3)\}$.
		\item $\Z_2\times\Z_2\times\Z_4$: $S_{11}=\{(0,0,1),\allowbreak(0,0,2),\allowbreak(0,0,3),\allowbreak(0,1,0),\allowbreak(1,0,0),\allowbreak(1,1,0)\}$, $S_{22}=\{(0,0,2),\allowbreak(0,1,0),\allowbreak(0,1,2),\allowbreak(1,0,1),\allowbreak(1,0,3),\allowbreak(1,1,1),\allowbreak(1,1,3)\}$, and $S_{12}=\{(0,0,0),\allowbreak(0,1,2),\allowbreak(1,0,3),\allowbreak(1,1,1)\}$.
		\item $\Z_2\times\Z_8$: $S_{11}=\{(0,1),\allowbreak(0,4),\allowbreak(0,7),\allowbreak(1,1),\allowbreak(1,7)\}$, \allowbreak $S_{22}=\{(0,1),\allowbreak(0,2),\allowbreak(0,3),\allowbreak(0,5),\allowbreak(0,6),\allowbreak(0,7)\}$, \allowbreak $S_{12}=\{(0,0),\allowbreak(0,2),\allowbreak(0,3),\allowbreak(0,4),\allowbreak(0,6),\allowbreak(0,7)\}$.
		\item $\Z_2\times\Z_8$: $S_{11}=\{(0,2),\allowbreak(0,4),\allowbreak(0,6),\allowbreak(1,2),\allowbreak(1,6)\}$, \allowbreak $S_{22}=\{(0,1),\allowbreak(0,2),\allowbreak(0,3),\allowbreak(0,5),\allowbreak(0,6),\allowbreak(0,7)\}$, \allowbreak $S_{12}=\{0\}\times\Z_8$.
		\item $\Z_2\times\Z_8$: $S_{11}=\{(0,2),\allowbreak(0,4),\allowbreak(0,6),\allowbreak(1,0),\allowbreak(1,2),\allowbreak(1,6)\}$, \allowbreak $S_{22}=\{0\}\times(\Z_8\setminus\{0\})$, \allowbreak $S_{12}=\{0\}\times\Z_8$.
		\item $\Z_4\times\Z_4$: $S_{11}=\{(0,1),\allowbreak(0,3),\allowbreak(1,0),\allowbreak(1,1),\allowbreak(3,0),\allowbreak(3,3)\}$, \allowbreak $S_{22}=\{(0,1),\allowbreak(0,2),\allowbreak(0,3),\allowbreak(2,0),\allowbreak(2,1),\allowbreak(2,2),\allowbreak(2,3)\}$, \allowbreak $S_{12}=\{0,2\}\times\Z_4$.
		\item $\Z_4\times\Z_4$: $S_{11}=\{(0,1),\allowbreak(0,2),\allowbreak(0,3),\allowbreak(1,0),\allowbreak(2,0),\allowbreak(3,0)\}$, \allowbreak $S_{22}=\{(0,2),\allowbreak(1,1),\allowbreak(1,3),\allowbreak(2,0),\allowbreak(2,2),\allowbreak(3,1),\allowbreak(3,3)\}$, \allowbreak $S_{12}=\{(0,0),\allowbreak(0,2),\allowbreak(1,1),\allowbreak(1,3),\allowbreak(2,0),\allowbreak(2,2),\allowbreak(3,1),\allowbreak(3,3)\}$.
		\item $\Z_{24}$: $S_{11}=\{\pm1,\pm5,\pm6,\pm7,\pm11\}$, \allowbreak $S_{22}=\{2,4,\dots,22\}$, \allowbreak $S_{12}=\{0,6,12,18\}$.
	\end{enumerate}
	
	Several of these graphs also arise over other groups: graph~1 over $\Z_2^3$ and $D_4$, graph~2 over $D_6$, graph~3 over $D_6$, graph~4 over $\Z_2^2\times\Z_3$ and $D_6$, graph~6 over $\Z_2^2\times\Z_4$, graph~12 over $\Z_4\times\Z_4$ and $\Z_2^2\times\Z_4$, and graphs~13 and~15 over $\Z_2^2\times\Z_4$.
	
	\paragraph{Further examples.} The graphs of Table~\ref{tab:new} that are not constructed in the main text are the following.
	\begin{itemize}[leftmargin=2em,itemsep=2pt]
		\item $n=24$: take the cube $Q_3$ and partition its vertices into the four antipodal pairs $P_0,\dots,P_3$. Add four sets $X_0,\dots,X_3$ of four new vertices each, join every vertex of $P_b$ to every vertex of $X_b$, and add all edges inside $X_0\cup X_3$ and inside $X_1\cup X_2$, so that these two sets induce copies of $K_8$. The cube vertices have degree $3+4=7$ and the new vertices $2+7=9$.
		\item $n=55$: take five cliques $K_5$, indexed by $b<5$, and join the clique with index $b$ completely to six new vertices $\ell_{b,1},\dots,\ell_{b,6}$; finally, join $\ell_{b,i}$ and $\ell_{b',i}$ whenever $b\ne b'$. The new vertices have degree $5+4=9$ and the clique vertices $4+6=10$.
		\item $n=60$: for $c\in\{0,1\}$ and $g\in\{0,1,2\}$ take an edge $L_{c,g}$ and a set $H_{c,g}$ of eight vertices. Each $H_c=H_{c,0}\cup H_{c,1}\cup H_{c,2}$ induces $K_{24}$, every vertex of $L_{0,0}\cup L_{0,1}\cup L_{0,2}$ is adjacent to every vertex of $L_{1,0}\cup L_{1,1}\cup L_{1,2}$, and $L_{c,g}$ is joined completely to $H_{c',g'}$ whenever $c\ne c'$ and $g\ne g'$. The vertices of the edges $L_{c,g}$ have degree $1+6+16=23$ and the others $23+4=27$.
		\item $n=192$: take three copies $C_0,C_1,C_2$ of $K_4$, completely joined to one another, and for each $i$ join $C_i$ completely to a copy $P_i$ of the complete tripartite graph $K_{8,8,8}$, which is in turn joined completely to a copy $D_i$ of $K_{36}$. The degrees are $3+8+24=35$ on $C_i$, $16+4+36=56$ on $P_i$ and $35+24=59$ on $D_i$. This is the blow-up of a triangle with a pendant path of length two at each vertex.
	\end{itemize}
	The graph on $33$ vertices of Table~\ref{tab:new} arose in a search over graphs with a two-orbit abelian group of automorphisms and is given in the ancillary files.
	
	\medskip
	\noindent\textbf{Data and reproducibility.}
	{\sloppy
		The accompanying supplementary archive is organised into three directories alongside the verification scripts. The \texttt{data} directory provides the \texttt{graph6} files \texttt{examples.g6} (the sixteen graphs of Table~\ref{tab:examples}), \texttt{examples\_new.g6} (the ten graphs of Table~\ref{tab:new}), \texttt{catalogue.g6} (all $138$ graphs of Section~\ref{sec:far}, with their tabulated invariants in \texttt{catalogue.tsv}), and \texttt{examples\_order10.g6} (the two ten-vertex graphs of Proposition~\ref{prop:tenvertices}), together with the parameter list \texttt{fourparam\_members.txt} for the $3654$ members of the family $\mathcal G(s,a,k,h)$ from Appendix~\ref{app:fourparam}. For each graph in Table~\ref{tab:new}, the \texttt{certificates} directory supplies an explicit integer matrix $X$ and integer $c$ satisfying $(nL+J)X=cI$, allowing resistance-regularity to be certified by a single integer matrix multiplication. Thirteen self-contained Python~3 and C scripts (\texttt{01}--\texttt{13}), using \texttt{numpy}, \texttt{networkx}, \texttt{sympy}, \texttt{python-flint}, \texttt{pynauty} and \texttt{nauty}/\texttt{geng} \cite{McKayPiperno2014}, verify all exact rational invariants and symbolic derivations across Sections~\ref{sec:forced}--\ref{sec:rows} and Appendix~\ref{app:fourparam}, including the built-in \texttt{geng} filter and per-chunk audit for $11\le n\le13$. Full run logs, command lines, software versions and the backward-error analysis are provided in \texttt{logs/} and \texttt{scripts/README.md}.\par}
	
	\setlength{\bibsep}{2.5pt plus 1pt}
	\bibliographystyle{ieeetr}
	\bibliography{non_regular_resistance_regular}
	
\end{document}